\documentclass{article}
\usepackage{amsmath,amsfonts,bm}

\usepackage{fontspec}

\usepackage{paralist}
\usepackage[left=2cm,right=2cm,top=2cm,bottom=3cm]{geometry}

\usepackage{natbib}
\usepackage[colorlinks=true,citecolor=blue]{hyperref}
\usepackage{amsthm}
\usepackage{amssymb}
\usepackage{physics}
\usepackage{booktabs}
\usepackage{comment}
\usepackage{tikz-cd}
\usepackage{tikz}
\usepackage{titlesec}
\titleformat{\chapter}[display]
  {\normalfont\huge\bfseries}
  {Appendix \thechapter:}{1em}{}
\usepackage{tikz-cd}
\usepackage{tikz}
\usepackage{pgfplots}
\pgfplotsset{compat=1.18}
\usepackage{scalefnt}
\usetikzlibrary{positioning}
\usepackage{subcaption}
\usepackage{tabularx}
\usepackage{stmaryrd,graphicx}
\usepackage{mathtools, mathrsfs}
\usepackage{braket}
\usepackage{algorithm}
\usepackage{algpseudocode} 
\usepackage[percent]{overpic}
\usepackage{subcaption}
\usepackage{xcolor}
\definecolor{ForestGreen}{rgb}{0.0, 0.3, 0.0}
\usepackage{makecell}

\newtheorem{theorem}{Theorem}[section]
\newtheorem{lemma}[theorem]{Lemma}

\theoremstyle{definition}
\newtheorem{definition}[theorem]{Definition}
\newtheorem{example}[theorem]{Example}
\newtheorem{hypothesis}[theorem]{Hypothesis}

\theoremstyle{remark}
\newtheorem{remark}[theorem]{Remark}
\newtheorem{corollary}[theorem]{Corollary}

\usepackage{thmtools}
\usepackage{etoolbox} 

\newcommand{\Q}{\mathbb{Q}}

\newcommand{\activation}{\phi}

\newcommand{\sumi}{\sum_{i=1}^N}

\newcommand{\CL}{\mathcal{L}}
\newcommand{\CX}{\mathcal{X}}
\newcommand{\CB}{\mathcal{B}}
\newcommand{\CY}{\mathcal{Y}}
\newcommand{\CN}{\mathcal{A}}

\newcommand{\CF}{\mathcal{F}}

\newcommand{\cC}{\mathcal{C}}
\newcommand{\cL}{\mathscr{L}}
\newcommand{\cP}{\mathscr{P}}

\newcommand{\R}{\mathbb{R}}

\newcommand{\N}{\mathbb{N}}
\newcommand{\eps}{\epsilon}
\newcommand{\weightvector}{{\mathbf{w}}}

\newcommand{\wzero}{\weightvector(0)}

\newcommand{\Ndw}{{\CN_{\weightvector}}}

\newcommand{\Ndwdag}{{\CN_{\weightvector}^\dagger}}

\newcommand*\FF{{\mathbf{F}}}

\newcommand{\FFg}{\FF[g]}

\newcommand*\GN{G_\infty}

\newcommand{\mrmD}{\mathrm{D}}
\newcommand*\mbbR{{\mathbb{R}}}

\newcommand*\NdwX{\CN^\CX_{\weightvector}}

\newcommand{\CNzero}{\CN(0)}

\newcommand{\CBep}{\CB_{\Phi}}
\newcommand{\CBMep}{\CB^M_{\Phi}}
\newcommand{\CBdw}{\CB_{\wst}}

\newcommand{\CNw}{\CN_{\weightvector}}
\newcommand{\CNwst}{\CN({\weightvector^*})}
\newcommand{\cNw}{\CN(\weightvector)}
\newcommand{\CNww}{\CN_{\weightvector\weightvector}}

\newcommand{\tmcX}{\widetilde{\CX}}

\newcommand{\Thetaw}{\Theta(\weightvector)}
\newcommand{\thetaw}{\theta(\weightvector)}

\newcommand{\Thetat}{\Theta(t)}

\newcommand{\eqdef}{\coloneqq}
\newcommand{\cLw}{\mathscr{L}[\weightvector]}
\newcommand{\DcLw}{\nabla \mathscr{L}[\weightvector]}

\newcommand{\CLN}{\CL[\cNw]}

\newcommand{\rmD}{\mathrm{D}}
\newcommand{\wst}{\weightvector^*}
\newcommand{\wsit}{\weightvector^*_{i}}
\newcommand{\wsiit}{\weightvector_{i+1}^*}

\let\epsilon\relax
\newcommand{\epsilon}{\varepsilon}

\newcommand{\Gi}{G_{M_i}}

\newcommand{\proj}{\text{proj}}

\newcommand{\Rinf}{\R^\infty}
\newcommand{\wt}{\weightvector(t)}
\newcommand{\CNt}{\CN(t)}

\newcommand{\wi}{\weightvector_{i}}

\newcommand{\Wkp}{W^{k,p}}

\providecommand{\keywords}[1]{\textbf{Keywords:} #1}

\begin{document}

\title{Man, Machine, and Mathematics}
\author{Akshunna S. Dogra\footnote{The NSF AI Institute for Artificial Intelligence and Fundamental Interactions (IAIFI), USA \\
$\textbf{ }\text{ }\quad$EPSRC Centre for Doctoral Training: Mathematics of Random Systems, UK\\
$\textbf{ }\text{ }\quad$Center of Mathematical Sciences and Applications, Harvard University, USA\\
$\textbf{ }\text{ }\quad$Department of Physics, Massachusetts Institute of Technology (MIT), USA \\
$\textbf{ }\text{ }\quad$Department of Mathematics, Imperial College London, UK\\
$\textbf{ }\text{ }\quad$MathePhysics, India. {\textbf{Correspondence}: adogra@nyu.edu}}
}
\maketitle
\begin{abstract}
Nonlinear models and optimization methods have successfully tackled a rapidly growing set of problems in recent years. Indeed, a relatively small toolbox of such models and methods can provide sufficient performance across a large landscape of tasks: deep learning alone has made significant recent contributions in scientific modelling, natural language processing, visual analysis, etc. A similar relationship exists between physical theories and phenomena, where many applications and observations emerge neatly from remarkably minimal foundations. It is natural to wonder if sparse unified frameworks could be built to steer discussion and discovery in the fields concerned with learning, optimization, and modelling. In this work, we posit and examine a possible outline for such a unified theory, interpreting the notion of “learning'' in a broad sense. In particular, we pursue our goals by viewing learning as an inter-connected process on multiple levels: problem setup, choosing methods, and the analysis of their interplay via imposed optimisation dynamics. We begin by proposing a precise yet versatile definition for “solvable'' problems. We then define the “parametrised methods'' by which their solution(s) may be “learned''. Our goal is to sketch a “universal convergence theorem'', specifying how and when solvable problems become amenable to the methods chosen for them. We find these constructions reduce the study of learning down to remarkably few ideas and tools - many of which are simply adapted from existing ones in dynamical systems theory, geometry, and fundamental physics.
\end{abstract}
\noindent \keywords{Mathematics of Learning; Physics of Learning; Gradient flows; Universal Approximation Theorem;$\quad$ Deep Learning; Smale-Hirsch Immersion theory; Lyapunov functions; Desingularization; Direct Limits.}

\section{Introduction}

The broad notion of “learning'' underpins biological activity across many levels of life. Indeed, in the human context, it acquires many-fold meaning: the biologically encoded notions of learned behavior, socio-economic activities, the machine based context added by modern computing, etc, all invoke it in myriad ways. We may see its signature in activities as disparate as mammals understanding spatio-temporal positioning, training a classifier, estimating functions using Neural/Kolmogorov Arnold Networks (NN/KAN), etc. Roughly speaking, it is the capacity of both machine and biological entities to incorporate information in their pursuit of some ill-to-precisely defined objective.

Within the computational context, a few modelling methods and heuristic principles are the cornerstones for a significant portion of modern optimisation tasks. Across computational disciplines as varied as natural language processing to numerically solving PDEs, the basic schema of finding \textbf{i)} a loss/control map, \textbf{ii)} some parameterization of the target space, and \textbf{iii)} an apt parameter update rule, is enough to generate sufficiently powerful models of desired solutions \cite{lecun15nat}. Indeed, many neural and biological mechanisms can also be viewed similarly \cite{dog25_2, freeenergy_brain_friston_07, control_motorcoord_todorov02}.

This bears strong resemblance to the influence of unified theories in physics - seemingly complex and unrelated phenomena emerging neatly from sparse and fundamental descriptions. Generalization and unification have been powerful forces in guiding the development of physics for centuries. Formal theories with a sparsity of rules describing a multitude of phenomena are considered objects of great beauty \cite{weinberg_dreams_93}, in addition to their obvious utilities vis a vis analysis and applications. This analogy motivates our search for a sparse and unified theory for learning.

Unsurprisingly, this goal is far too ambitious and ill-defined. We will substantially narrow our focus to characterize learning pursued via “parameterized methods'', where a model and its parameters evolve under some optimisation rule to produce better estimates of the desired target/solution. Even this goal is too broad and ill-defined, but we can now pose questions with a hope of being answered. The following ones lie at the heart of our pursuit:
\begin{enumerate}
    \item To what extent can we formalise the notion of a solvable problem and a parametrised method? (Sec. \ref{prob_prec}, \ref{arch_prec})

    \item How do we use parametrised methods to solve problems and why are nonlinear ones useful? (Sec. \ref{Optimization_Flows}, \ref{nlin_A})

    \item How do we make a nominally solvable problem more tractable to a chosen method? (Sec. \ref{FNtogether}, \ref{arch_exp_prun})
\end{enumerate}

The moonshot result, inspired from universal approximation theorems prevalent in machine learning \cite{Horn91, UAP_Mishra_21}, would be a “universal convergence theorem'' for learning. It would specify exactly how/when a problem solvable \textit{in principle}, is tractable to the chosen method \textit{in practice}. The goal of this work is to sketch what prerequisites are needed for that result and what it might look like: outline what we want from a general theory of learning and why. 

We have broken the notion of learning into three cornerstones that also supply our title (Fig. \ref{MMMTriangleLearning}\textbf{a}): the problem (Man), the method we choose to attack it with (Machine), and the formal analysis of their interplay (Mathematics). Yet, somewhat poetically, a lot of the connective tissue between these facets is also supplied by ideas originating in physics: the tools needed for learning about our universe often being the tools that describe how we learn.

Gradient Flows and their variants/generalizations are one of the linchpins of modern optimisation and modelling paradigms. On an analogy level, they are to learning as Schrödinger's equation is to Quantum Mechanics and for reasons beyond just their ubiquity (Remark \ref{loss_para}: most forms of \textit{good} optimisation are \textit{equivalent} to some gradient flow). We present a brief outline of the technical challenges and aims by scaffolding the discussion around them.

\begin{figure}
\centering

\begin{minipage}{0.3\textwidth}
\centering
\begin{tikzpicture}[
    every node/.style={draw, circle, minimum size=1.2cm},
    >=Stealth,
    thick
]
\node (A) at (90:2)  {{Problem}};
\node (B) at (210:2) {{$\,$Method$\,$}};
\node (C) at (330:2) {{Analysis}};
\draw[<->, bend right=30] (A) to (B);
\draw[<->, bend right=30] (B) to (C);
\draw[<->, bend right=30] (C) to (A);
\node[draw=none, shape=rectangle] at (0,0) {\textbf{Learning}};
\end{tikzpicture}
\end{minipage}
\hfill
\begin{minipage}{0.67\textwidth}
\centering
\renewcommand{\arraystretch}{1.1}
\begin{tabular}{|l|l|l|}
\hline
\textbf{Problem of interest} & \textbf{Common methods} & \textbf{Fields of analysis} \\
\hline
Regression Analysis  & Linear regression & Linear Algebra  \\
\hline
Function Estimation  & Fourier Series, KANs & Fourier Analysis  \\
\hline
Spatial Navigation & Mammalian grid cells & \makecell[l]{Signal processing,\\ Neuroscience} \\
\hline
\makecell[l]{Shape and Visual\\Recognition } & \makecell[l]{Wasserstein flows with\\parametrised measures} & \makecell[l]{Computer Vision,\\ Optimal Transport} \\
\hline
Classification  & Support Vector Machine & Manifold learning \\
\hline
Differential Equations & Finite Element Methods & Numerical DEs \\
\hline
Text Generation  & NNs, Transformers & NLP, Linguistics  \\
\hline
Human Communication & Several neural pathways & Neuroscience \\
\hline
\end{tabular}
\end{minipage}

\caption{(\textbf{a}) Schema for learning, \hfill (\textbf{b}) Some common problems,  usual methods, and invoked fields of analysis.$\,\,\,$ }
\label{MMMTriangleLearning}
\end{figure}

\subsection{An informal summary of the desired results, motivations, and layout of this work}\label{informal_sum}

\noindent Our first task is generating formal proxies for ``well-behaved'' or ``solvable'' problems. Roughly speaking, problems are represented by maps $\FF: G \to H$ acting between apt metric spaces. We call them well-behaved if there exists a ``solution(s)'' $\Phi \in G$, a neighbourhood $\CBep \subset G$, and a loss functional $\CL: G \to \mathbb{R}$, such that something resembling the following ``gradient flow'' can be well-defined,  converges to $\Phi$, and minimises both $\CL$ and $|\FF|$ (Defn. \ref{well-behaved}):
\begin{equation}\label{rough_gradientFlow}
\frac{d}{dt}g(t) = -\nabla \CL[g], \qquad\qquad \qquad\qquad \CBep \ni g(s)_{s \in \R^+}  \implies g(t) \xrightarrow[t \to \infty]{} \Phi
\end{equation}
Assuming $G$ has sufficient structure to support differentials, variants of the Lojasiewicz inequality\footref{General_LI} provide some of the weakest known conditions on $\CL$ and Eq. \ref{rough_gradientFlow} for such convergence guarantees. We note its simplest form below:
\begin{equation}\label{informal_LI}
    |\CL[g] - \CL[\Phi]|^\alpha \leq C\|\nabla \CL[g]\|, \qquad g \in \CBep,\textbf{ } \alpha \in [1/2, 1), \textbf{ }C > 0,
    \tag{LI}
\end{equation}
However, even if \ref{informal_LI} holds, often $\dim(G) = \infty$ and Eq. \ref{rough_gradientFlow} is computationally inaccessible. We need well-defined methods for generating parametrised models in $G$ and leveraging finite dimensional gradient flows over them. Architectures, $\CN: \R^M \to G$, use $M$ parameters to produce such model sets $G_M = \{\cNw: \weightvector \in \R^M\} \subset G$. If we set $\cL \eqdef \CL \circ \CN$ (a natural and fundamental choice - Remark \ref{loss_para}), we get some other gradient flows of interest:
\begin{subequations}\label{Rough_parameteric_gradient_flow}
\begin{alignat}{2}
    &\,\,\,\quad\frac{d}{dt}\wt = -\DcLw, \qquad\qquad\qquad\qquad\qquad\quad\qquad\quad \wt \in \R^M, \\
    \implies& \frac{d}{dt}\CN(\wt) = -\Ndw\Ndwdag \nabla \CL[\CN(\wt)], \qquad\qquad\qquad\qquad \CN(\wt) \in G_M
\end{alignat}
\end{subequations}
In general, a method is comprised of an architecture $\CN$ parametrising $G$, its parameters $\weightvector$, and the update rules for $\weightvector$. However, \textbf{good} update rules for minimising $\CL$ are \textbf{equivalent} to the system above \cite{Barta12}, so it will be our focus.

If $\CN[\wt]$ is eventually in $\CBep$, the architecture is ``well-initialised'' at $\wzero$ (Defn. \ref{well_initialize}). Broadly speaking, if $\FF$ is well-behaved and $\CN$ is linear and well-initialised, then it is easy to show that $\cL$ is at least ``locally'' well-behaved and the optimization dynamics are simple. In contrast, for nonlinear $\CN$, the optimization dynamics can be significantly more complex. Unfortunately, linearity is limiting  and produces at most $M$ dimensional linear subspaces in $G$.

This presents a significant issue if $\dim(G) > M$ and/or $\Phi$ is unknown, since we can't force the distance $d(\Phi, G_M)$ to be small. In contrast, nonlinear $\CN$ can allow us to ``cheat'' using the same number of parameters by ``folding'' $G_M$ in $G$ more efficiently. For example, consider how well a spiral covers $\R^2$ in comparison to a line, even though both are smooth 1-dimensional curves (Fig. \ref{fig:Lin_vs_NLin}). However, nonlinear $\CN$ may lead to $\cL$ that \textbf{(i)} have too many critical points, \textbf{(ii)} do not trivially inherit being well-behaved from $\CL$, and/or \textbf{(iii)} permit non-limiting dynamics.

Thm. \ref{vtk_central} addresses these issues for many methods: it establishes versatile conditions on $\CN$ to answer if/when well-behaved problems are also locally well-behaved with respect to (w.r.t.) the parameters of $\CN$. Specifically, it establishes the kinds of $\CN$ for which the critical points $\{\wst: \nabla\cL[\wst] = 0 \}$ of $\cL$ inherit a modified \ref{informal_LI} from $\CL$:
\begin{equation}
\exists \textbf{ } \alpha^* \in [1/2, 1), \textbf{ }C^* > 0, \textbf{ }\CBdw \subset \R^M,   \quad \text{ s.t. } \quad  |\cLw - \cL[\wst]|^{\alpha^*} \leq C^*\|\DcLw\|, \quad \weightvector \in \CBdw,     
\end{equation}
guaranteeing convergent parametric and model gradient flows $\wt \to \wst \in \R^M$ and $\CNt \to \CNwst \in G_M$. 

However, even with Thm. \ref{vtk_central}, if $M < \infty$ and $\dim(G) = \infty$, there is no guarantee we can place $\Phi \in G_M$ in a way that lets $\CN(\wt) \to \Phi$, especially as nonlinear $\CN$ can generate an excess of critical points for $\cL$. Fortunately, if $\CN$ is well-initialised, Thm. \ref{vtk_central_2} gives a chained, computationally feasible path to $\Phi$. It does so constructively, by gradually and iteratively expanding $\CN, G_M$ by adding new parameters {\textbf{without}} changing the model $\CNwst$. While a bit limited due its reliance on differentiability, Thm. \ref{vtk_central_2} is a good prototype of a universal convergence theorem.

We may summarise the above as a meta-heuristic, which we posit covers a large universe of learning processes:
\begin{enumerate}
    \item \textbf{Problem Setup}: Finding maps $\FF: G \to H, \CL: G \to \R^+$ between apt spaces that generate the setting in which to discuss the problem and implicitly or explicitly define a solution(s) $\Phi \in G$ w.r.t. $\FF, \CL$ (Sec. \ref{Tablesetting}).
    \item \textbf{Method Setup}: Picking an architecture $\CN: \R^M \to G$ s.t. the set of models $G_M$ is sufficiently dense and efficient in covering $G$, constructing an apt parametric form for $\CL$, and finding a good $\weightvector$ update rule (Sec. \ref{architectures}).
    \item \textbf{Mathematical Analysis}: Showing the losses $\CL$ and $\cL$ are well-behaved enough to allow useful optimisation flows, often by showing $\FF$ is well behaved and extending it to $\cL$ via properties of $\CN$ (Thm. \ref{vtk_central}).
    \item \textbf{Optimization/Error Correction}: Implementing constraint and analysis guided parametric optimisation and boosting existing models and methods (Thm. \ref{vtk_central_2}).
\end{enumerate}

Our theme is learning and our focus is on solvable problems and parametrised methods: naturally, scientific machine learning (ML) and neuroscience are strong inspirations for us. However, note the seemingly absent role of ``data'' above - it  will be our critical departure from conventional ML analysis. Data plays a critical role in determining how a problem, method, etc are setup. However, while its influence is implicitly strong in every choice and computation, it is not an explicit variable of the framework post those choices.

Our focus will be pedagogy and conjecture, \textbf{not} proof making or numerical validation. However, we rely heavily on and borrow liberally from a rich literature: the references, conclusions, and appendix will be used to give a full accounting of the invoked concepts, supporting empirical results, proofs, etc (in particular, App. \ref{nomenclature} summarises notation, App. \ref{pert_note} contains supplemental comments, and App. \ref{fund_lemmas} details all key concepts in full rigour). 

This work, derived from a series of lecture notes, is equal parts a survey, scientific manifesto, and a set of results and conjectures. Superficially, the universe of problems, methods, and invoked fields of analysis is vast (Fig. \ref{MMMTriangleLearning}\textbf{b}). We are claiming there are  common threads running through them all. Let us begin by defining ``solvable'' problems.

\section{The Problem (Man)}\label{Tablesetting}

\begin{figure}[h]
\centering

\begin{minipage}{0.3\textwidth}
\centering
\begin{tikzpicture}[
    scale=1.2,
    transform shape,
    every node/.style={
        draw, 
        circle, 
        minimum size=1cm, 
        inner sep=2pt,    
        align=center       
    },
    >=Stealth,
    thick
]
\node (A) at (90:2) {$\,\,$Target$\,\,$\\ $\,\,$Spaces$\,\,$};
\node (B) at (210:2) {Nominal\\Loss};
\node (C) at (330:2) {$\,$Solutions};

\draw[<->, bend right=30] (A) to (B);
\draw[<->, bend right=30] (B) to (C);
\draw[<->, bend right=30] (C) to (A);

\node[draw=none] at (0,0) {\textbf{Problem}};
\end{tikzpicture}
\end{minipage}
\hspace{75pt}
\begin{minipage}{0.479\textwidth}
\centering
\begin{overpic}[width=\textwidth]{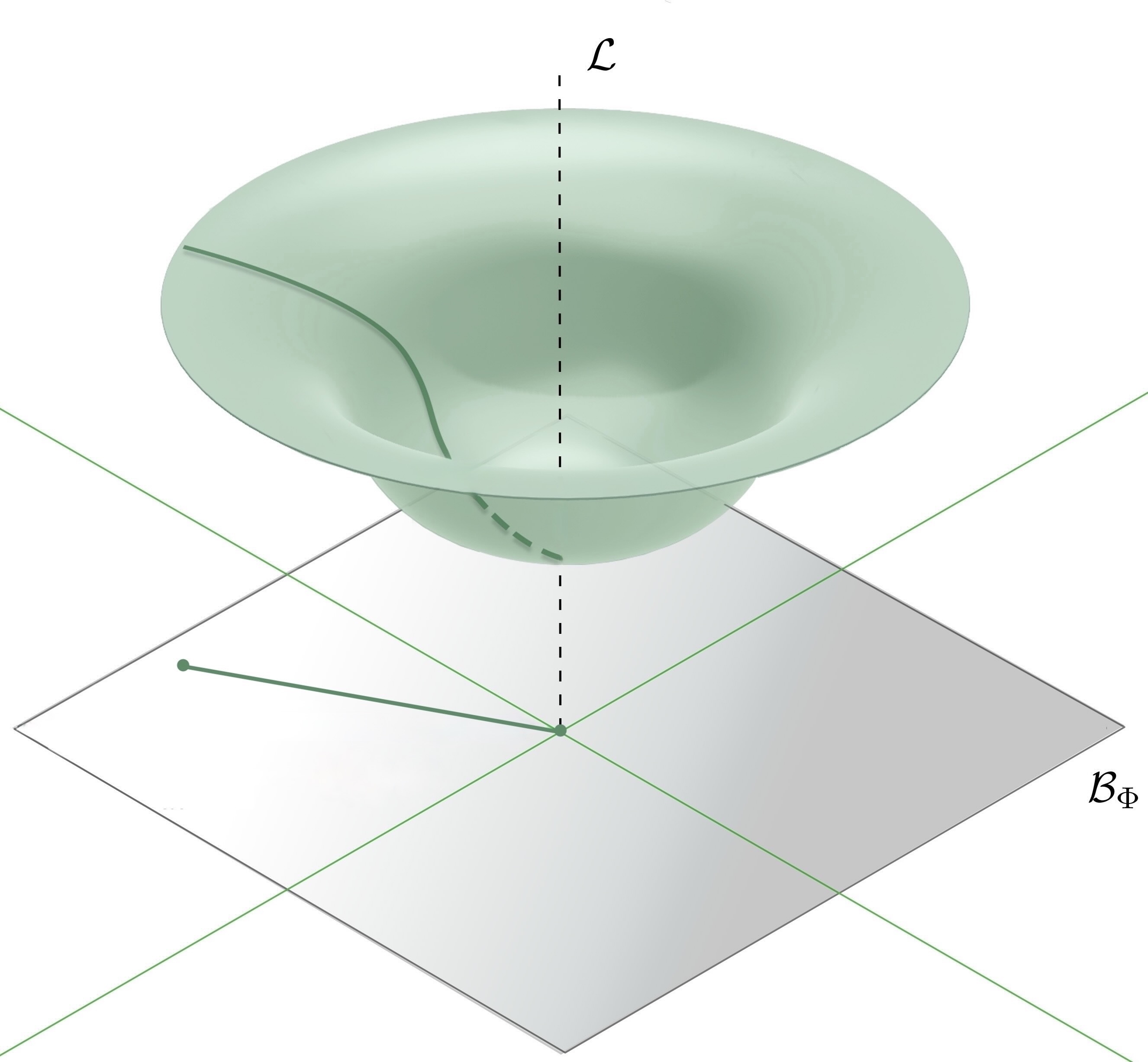}
        \put(47.5,24.5){\makebox[0pt][l]{{$\Phi$}}}
        \put(20.2,35){\makebox[0pt][l]{ {${g(t)}$}}}
        \put(53,28){\makebox[0pt][l]{{${\frac{d}{dt}g(t) \, = \,- \nabla\CL[g(t)]}$}}}
        \put(64, 34){\makebox[0pt][s]{{${g(0) \in \CBep}$}}}
\end{overpic}
\end{minipage}

\caption{(\textbf{a}) A generic schema for problems, \hfill (\textbf{b}) Nominal gradient flow dynamics for a solvable problem.}
\label{MMMTriangle}
\end{figure}

The notion of “learning'' is intimately connected to the notion of having a “problem'' to be solved or learned upon. We may wish to classify objects \cite{trajan25}, understand visual cues \cite{dog23, feydy_20}, produce solution models for differential equations \cite{dog23_3, weinan18, amina21, mar20, shin_karni_20, Sirig18}, chat with a human being \cite{Attention_all_need_Ilya_17}, etc. Thus, if we are to tackle “learning'', we must first precisely define what kinds of problems we are willing to consider and what it means for them to be solvable.

We restrict ourselves to the settings where to have learned something is to have identified something “correct'' from a list of possibilities - to have picked out a solution $\Phi$ from a target space $G$. The target space could be a Euclidean space, a space of measures, a space of smooth maps between other input/output spaces $\CX, \CY$, etc. In turn, the fitness of our choice has to be measured by some proxy - say a loss function $\CL$ - defined over this target space. In essence, $\CL$ tells us how good a particular $g \in G$ is relative to our desired solution $\Phi$. Solvability is then simply a condition on the dynamics permitted on $G$: on this level, dynamics leading to $\Phi$ is the process of learning. 

We need to convert this line of reasoning into a versatile yet precise mathematical proxy (as sketched in Sec. \ref{informal_sum}). However, before getting into specifics, is this a sufficiently powerful structure to support many kinds of “learning'' (computational and otherwise)? As an example, when solving differential equations, this is often how existence and well-posedness is shown: one defines an appropriate variational energy $
\CL \in \cC^2(G)$ over some space of functions $G$, shows that a unique global minima exists, and that $\nabla \CL$ is invertible around (if not at) it \cite[Ch. 8.2]{evans1998}\cite{manyfold_learning, mar20}. 

The above example hints that to build a unified theory, we need the weakest convergence conditions on gradient flows. Further, if $\CL$ is an integral associated with $\FF$ (as variational energies are), then we may be able to define the derivative $D \CL$ as well. Then all we need is some natural isomorphism between $D\CL$ and the desired “gradient'' element $\nabla \CL$. This points us to what $G$ has to be. But is this differentiable loss + convergent flow pairing versatile?

In \cite{manyfold_learning, feydy_20}, shape and visual recognition problems are interpreted as a problem of matching model measures with target measures. Signals over some detector geometry (pixel intensities on a screen or energy depositions on CERN detector surfaces) are interpreted as measures. The set of these measures can be given a manifold-like geometry using Wasserstein metrics \cite{Ott01,Lot08}. In turn, this directly led \cite{dog23, manyfold_learning, feydy_20} to loss functions and gradient flows to match the model and target measures. Classification problems can be given similar structure using Bochner spaces \cite{trajan25}.

In neuroscience, control theory provides a similar and fruitful unifying perspective on how to view many neuron modulated actions \cite{freeenergy_brain_friston_07, control_motorcoord_todorov02}. Similarly, mammalian navigation is modulated by ensembles of grid cells that fire according to Gaussian rules based on the current location and the instinct to go to a privileged location, measured by how well the correlation function at the current location matches the correlation at the desired location. In \cite{dog25_2}, it was shown the issue of identifying the correct location can also be reduced to searching for the extremum of a correlation function, which unsurprisingly is best done through an analysis of its derivative w.r.t. the location.

Transformers are simply sequence-to-sequence mappings over token representations \cite{Attention_all_need_Ilya_17}: piecewise differentiable maps between high but finite dimensional Euclidean spaces. Even on feature level analysis, the optimization of transformers can be studied by invoking Wasserstein gradient flows \cite{Wass_training_transformers_Gao_24}. We may even stretch these arguments to cover human communication, under the simplification it is reception/reaction to audio-visual cues.

Regardless, differentiable losses and gradient flows over differentiable manifolds are clearly a common element across many domains (see \cite[Ch. 3]{manyfold_learning} for complete discussions and proofs for many of the arguments made above). Further, solvability can be tied to the existence of gradient flows that are ``well-behaved'' in the neighbourhood of ``solution(s)''. We make these ideas concrete using the Lojasiewicz Inequalities \cite{isobe_LI_Wass_23, Loja_og_63, Loja_og_65, Loja_og_83} as a hypothesis on $\FF$ and $\CL$, which ties the geometry of our problems to the dynamics possible on their target spaces.

\begin{table}
\centering
\renewcommand{\arraystretch}{1.25}
\begin{tabular}{|l|l|l|l|}
\hline
\textbf{Problem} 
& \textbf{Typical Target Space $G$} 
& \textbf{Typical Loss} 
& \textbf{Typical solvability} \\
\hline
Regression Analysis \cite{DraperSmithRegression_98} & Space of smooth maps:  $\cC^\infty$ & $\CL[g] = \braket{g - \Phi, g - \Phi}$ & Strongly convex $\CL$ \\
\hline
Function Estimation \cite{SteinShakarchiFourier_03} & Lebesgue spaces: $L^p(\CX)$ & $\CL[g] = \braket{g - \Phi, g - \Phi}$ & Strongly convex $\CL$ \\
\hline
Spatial Navigation \cite{dog25_2} & Euclidean spaces: $\R^2$ or $\R^3$ & Correlation distances & Lojasiewicz Inequality \\
\hline
\makecell[l]{Shape and Visual\\ Recognition \cite{manyfold_learning, feydy_20}} & \makecell[l]{Smooth Wasserstein spaces\\ \cite{Lot08}: $\cP^\infty_2(\CX)$ or $\cP^\infty_1(\CX)$} & \makecell[l]{2-Wasserstein metric,\\1-Wasserstein metric} & \makecell[l]{Wasserstein LI \cite{LI_Wasserstein_18},\\ Strongly convex $\CL$} \\
\hline
Classification \cite{trajan25, KL_SCHWAB_06} & Bochner spaces: $L^p(\CX; \CY)$ & Bochner norm &  Strongly convex $\CL$ \\
\hline
Differential Equations \cite{dog23_3} & Sobolev spaces: $W^{k, p}(\CX)$ & $\CL = \braket{\FFg, \FFg}_H$ & LI, Strongly convex $\CL$ \\
\hline
\makecell[l]{Text Generation \cite{Attention_all_need_Ilya_17}\\ (Distributional view \cite{Wass_training_transformers_Gao_24})}& \makecell[l]{seq2seq space: $\cC^1(\R^{dn}, \R^{dm})$\\(Smooth Wasserstein spaces)} & \makecell[l]{Cross Entropy, others\\(Weight decay losses)} & LI, Strongly convex $\CL$ \\
\hline
\end{tabular}
\caption{Problems, commonly used loss functions, and their solvability conditions}
\label{tab:li_summary}
\end{table}

\subsection{Solvable or Well-behaved Problems}\label{prob_prec}

\begin{hypothesis}[Well-behaved]\label{well-behaved}
Assume $\FF:G \to H$, where $G$ is a complete, connected, separable, $\cC^2$-Banach-Riemannian manifold\footref{Banach_manifold} and $H$ is a Hilbert space. $\FF$ is said to be \textbf{well-behaved} w.r.t. a nominal loss functional $\CL:G \to \R^+$ and a solution $\Phi \in G$, if \textbf{(i)} $\braket{\FF[g], \FF[g]}_H$ and $\CL[g]$ are both minimized at $g = \Phi$, \textbf{(ii)} $\CL$ is coercive\footref{coercive}, \textbf{(iii)} there exists a neighbourhood $\CBep \subset G$ s.t. $\CL \in \cC^2(\CBep, \R^+)$, and \textbf{(iv)} $\CL$ satisfies the following gradient \ref{Lojasiewiczineqaulity}:
\begin{equation}\label{Lojasiewiczineqaulity}
\exists \alpha \in [1/2, 1), \textbf{ }C > 0, \text{ s.t. } \forall g \in \CBep, \qquad \qquad |\CL[g] - \CL[\Phi]|^\alpha \leq C\|\nabla \CL[g]\|_{T_g G}
\tag{LI}
\end{equation}
Further, $\FF$ is \textbf{strictly/strongly} well-behaved w.r.t. $\CL$, if $\CL$ is strictly/strongly convex\footref{strong_convex} over $\CBep$. Finally, $\FF$ is \textbf{locally} well-behaved w.r.t. $\CL$ at some $\Phi \in G$, if \textbf{(i)} - \textbf{(iv)} hold locally on $\CBep$.
\end{hypothesis}

\begin{remark}[Nominal Loss]\label{nominal_loss}
$\CL$ is called the “nominal'' loss as computers can't directly support gradient flows over $G$ if $\dim(G) = \infty$. $\CL[g] \eqdef \frac{1}{2}\braket{\FF[g], \FF[g]}_H$ is a common (not necessarily predominant) choice for nominal losses. Similarly, $G$ will usually be a Hilbert space. To ease discussion, we will also assume w.l.o.g. that $\CL[\Phi] = 0$.
\end{remark}

Each requirement serves a key purpose. Computers can only handle finite operations: $G$ is taken to be separable\footref{UAP_sep} so it is rich enough to allow generalisations to be made from countable operations over it. It is given Riemannian structure so it can support gradients. $G$ is often a Hilbert space, which is useful since compact operators over such $G$ have the approximation property, which allows numerical estimation of objects we will soon find key to modelling. However, Hypothesis \ref{well-behaved} concerns problems in the weaker settings as well (SHAPER \cite{dog23, manyfold_learning, feydy_20}).

Regularity of $\CL$ supports the gradient flows that could lead to $\Phi$, while coercivity\footref{coercive} ensures that $\CL[g(t)]$ limiting to its global minimum implies $g(t) \to \Phi$. The gradient \ref{Lojasiewiczineqaulity} (and its equivalent distance form\footref{grad_like_flow}) guarantees eventual convergence of gradient-like flows to $\Phi$ with \textit{a priori} or \textit{in situ} predictable rates if we start near $\Phi$ - even strict convexity can lack such convergence guarantees. Further, convexity\footref{vanilla_con} requires that we always get closer to $\Phi$ when $\CL$ decreases, while \ref{Lojasiewiczineqaulity} can allow temporary moments where $\CL$ decreases while we get \textit{further} from $\Phi$. 

Note our claim is not that all (or even most) problems are best solved under the framework established by Hypothesis \ref{well-behaved}: simply that most solvable ones can be setup in a form that satisfies it (at least nominally). \cite[Ch. 3]{manyfold_learning} rigorously establishes how, and under which $\CL$, most of Table \ref{tab:li_summary} satisfies Hypothesis \ref{well-behaved}. The justification for spatial navigation can be inferred directly from \cite{dog25_2}; \cite{isobe_LI_Wass_23} supplies it for the distributional view on text generation. 

Our running example for discussions in this work is a simplified nonlinear Poisson-Boltzmann equation (nPBE). 

\begin{example}[nPBE for a fixed potential $h \in L^{2}(\text{[}-\pi,\pi\text{]}^3)$ with zero boundary conditions]\label{stock}
\begin{equation}
\FF : W^{2, 2}([-\pi,\pi]^3) \cap W^{1, 2}_0([-\pi,\pi]^3)  \to L^{2}([-\pi,\pi]^3), \quad \FF[g] = -\Delta [g] + \sinh(g) + h
\end{equation}
\end{example}
This nonlinear PDE is key to problems in bio-medicine and chemistry \cite{iglesias2022weak}. It is defined between $\infty$-dimensional spaces and showcases intuitive features like the strong convexity\footref{strong_convex} of $\CL[g] = \braket{\FF[g], \FF[g]}_H$ near $g = \Phi$ (\cite[Ch. 3.2.1]{manyfold_learning}). It is regular enough to support gradient flows, so a good initial model $g$ could be evolved to the solution. Finally, it is mathematically tractable enough to lower bound the exact optimization rates we should expect.

The nominal losses $\CL[g]$ allow us to define an analytically instructive gradient flow:
\begin{equation}\label{idealized_G_flow}
\frac{d}{dt}g(t) = -\nabla \CL[g(t)],  \qquad\qquad\qquad g(t) \in G 
\end{equation}
For Example \ref{stock}, we have $
\nabla \CL[g] 
= 2\, (I - \Delta)^{-2} ( -\Delta + \cosh(g) ) (-\Delta g + \sinh(g) + h )$. Note that the gradient here is computed w.r.t. the $W^{2, 2}$ norm rather than the $L^2$ norm because we want convergence in $G$, not in $H$. Such choices are significant: $\nabla_{L^2} \CL[g] = ( -\Delta + \cosh(g) ) (-\Delta g + \sinh(g) + h )$, which would force us to pick $G$ with stronger regularity. “Analysis'' is critical even while we are setting up: the three cornerstones of Fig. \ref{MMMTriangleLearning} inform each other.

In theory, well-behaved problems are nice and optimize at determinable rates under Eq. \ref{idealized_G_flow}:
\begin{lemma}[\cite{Loja_Chill_convg_06}, Remark 2.8]\label{simple_GD_G}
Assume $\FF$ is well-behaved w.r.t. $\CL$ with \ref{Lojasiewiczineqaulity} constants $C, \alpha$ and $g(0) \in \CBep$. Under Eq. \ref{idealized_G_flow}, we have:
\begin{equation}\label{LI_convg_rates_basic}
g(t) \xrightarrow[t \to \infty]{} \Phi,    \qquad \|g(t) - \Phi\| = \begin{cases}
                O(e^{-Ct}), \text{ if } \alpha = \frac{1}{2}\\
                O(t^{-\frac{1 - \alpha}{2\alpha - 1}}), \text{ if }  \frac{1}{2} < \alpha < 1\\
                \end{cases}
\end{equation}
\end{lemma}

If we have additional information like strong convexity, our task is even easier. For Example \ref{stock}, we can show that $\alpha=1/2, C \geq 1$ for all minima of $\CL$ over $G$ \cite[App. C.1]{manyfold_learning}, which imbues strong convergence properties on the nPBE. For most non-PDE related applications and nominal losses, such results are even easier to get. Indeed, Lemma \ref{simple_GD_G} is a significantly weaker result than we can obtain \cite{kurdyka98}. Unsurprisingly, it is often possible (and even trivial) to establish that a problem of practical interest is solvable \textit{in principle}. But can Eq. \ref{idealized_G_flow} solve it \textit{in practice}?

Unfortunately, Eq. \ref{idealized_G_flow} can't be directly realized on a computer since they cannot handle generic infinite dimensional dynamics. We need a computationally accessible parameterization of $\CL$ and $G$ that inherits the convergence properties of Eq. \ref{idealized_G_flow}, even when using parameters nonlinearly. However, to make sense of that statement, we need to define parameterization, architectures, update rules, and the trade-off between linear vs nonlinear methods.

\section{The Method (Machine)}\label{architectures}

Now that we have established a definition for “solvable problems'', let us consider the methods used to solve them. Our focus will be on parametrised methods: those seeking to estimate a solution by iteratively refining some countable number of parameters to produce better and better models. We borrow AI/ML terminology to call the mapping that parametrises $G$ using $M$ parameters, an “architecture'' $\CN: \R^M \to G$. Its image $G_M = \{\cNw: \weightvector \in \R^M \}$ is called the model set for obvious reasons. A “method'' is then simply characterised by the chosen architecture(s) $\CN$, its parameters $\weightvector \in \R^M$, its models $\cNw \in G_M$, and the rules for evolving $\weightvector$ (Fig. \ref{fig:method_and_summary}\textbf{a}).

Classical modelling methods are often linear w.r.t. parameters, even when comprised of nonlinear expressions: indeed, many such methods like regression, Fourier analysis, etc, comprise their own well-established fields. However, both ML and biological learning leave this paradigm by using nonlinear parametrisations (Fig. \ref{fig:method_and_summary}\textbf{b}), usually through a base function (or activation) $\phi$. We will need results that exemplify the balance between these two choices.

The fitness of a method can be judged by how densely and efficiently its chosen model sets cover $G$. Consider $d(G_M, \Phi)$, the distance between the best $M$ parameter model of $\CN$ and $\Phi$. As $M \to \infty$ for a class of choices, we want $d(G_M, \Phi) \to 0$. Fortunately, universal approximation theorems already provide this in many\footref{UAP} settings \cite{GKP20, UAP_Wasserstein_space_Lu_20, UAP_Mishra_21}. Further, it is relatively simple to ensure that the model set of a linear $\CN$ is at least $M$ dimensional within $G$ ($M$-rank matrices are dense within the space of $M \times M$ matrices). As such, a nonlinear method is more useful if and only if $G_M$ has dimension $M$ (at least mostly), \textbf{while} offering denser coverage than a linear $\CN$ does. 

Euclidean spaces $\R^M$ will serve as the parameter spaces of choice when $M < \infty$. When $M = \infty$, two choices of considerable interest emerge: $l^2$ (the space of square summable sequences) and $\Rinf$ (the inductive limit of convergence\footref{inductive_limit} on $\R^M$ - the set of all sequences that are \textit{eventually} non-zero for only a finite number of their elements). However, as $\Rinf \subset l^2$ and $l^2$ is a Hilbert space, the reader may keep $l^2$ in mind for most of what follows.

Once an architecture $\CN$ is chosen, the next question is how we are to refine its parameters. The natural choice is to use $\weightvector$ update rules that expect to lower the value of $\CL$ every time $\weightvector$ is refined. The properties of $\cLw \eqdef \CL\circ \cNw$ and its gradient w.r.t. $\weightvector$ are critical in such considerations, even if we do not use them explicitly. 

Note again that we are making a strong distinction between \textit{parameters} and \textit{data}: the parameters produce models: these models may or may not then be refined using available data. $\R^M$ represents the parameter space, while $\CX \subseteq \R^d$ will represent the input data space if it is relevant. As a practical matter, the choice of $M$ and the parametric loss functions may or may not be related to $d$ or the number of available datapoints (say $N$). However, our focus is on the parametric/model dynamics once the choices are made - the role of data is implicit, not explicit.

In summary, the definition of an architecture should incorporate well studied linear (or near-linear) techniques, but with room for the wide variety of modern nonlinear methods. It should be well-suited to working with a finite collection of tunable parameters (a constraint of working with computers), but capable of extension to a countable number of parameters (so we can take limits). It should be regular enough to support gradient flows, but rich enough to produce model sets that cover complete separable spaces with sufficient density, if given enough parameters.

\begin{figure}
\centering

\begin{minipage}{0.25\textwidth}
\centering
\begin{tikzpicture}[
    every node/.style={draw, circle, minimum size=1.6cm, align=center},
    >=Stealth,
    thick
]
\node (A) at (90:1.7)  {Param-\\eters};
\node (B) at (210:1.7) {Model \\ sets };
\node (C) at (330:1.7) {Update\\rules};

\draw[<->, bend right=30] (A) to (B);
\draw[<->, bend right=30] (B) to (C);
\draw[<->, bend right=30] (C) to (A);

\node[draw=none, shape=rectangle] at (0,0) {\textbf{Method}};
\end{tikzpicture}
\end{minipage}
\hfill
\begin{minipage}{0.687\textwidth}
\centering
\renewcommand{\arraystretch}{1.25}
\small
\begin{tabular}{|l|c|c|l|}
\hline
\multicolumn{1}{|c|}{\textbf{Method}} 
& \multicolumn{2}{c|}{\textbf{Linear in}} 
& \multicolumn{1}{c|}{\textbf{Typical Models $\cNw$}} \\
\cline{2-3}
&  $\weightvector$ &  $(\cdot)$ & $N$: number of data points, $d$: data dimension \\
\hline
Linear Regression 
& Yes & Yes 
& $ (w_{ij})_{i, j \in \{1, 2,..., \sqrt{M}\}}$; (Matrices) \\
\hline
Fourier Series 
& Yes & No
& $ \sum_{m=-\infty}^{\infty} w_m e^{i2\pi \frac{m}{P}(\cdot)}$; (Functions, $M = \infty$) \\
\hline
KANs \cite{KAN_25}
& No & No
& $ \sum_{i=1}^{2d} \psi_i (\sum_{j=1}^d \sum_{k = 1}^{W}  w_{ijk}\phi_{jk}(x_j)); M \sim  Wd^2$ \\
\hline
Grid-cell Neurons 
& No & No 
& $ \frac{1}{N}\sumi \sum_{j\in \N} \exp{-(\frac{\weightvector - j \lambda_i}{2c_i})^2}; \weightvector \in \R^M = \R^d$ \\
\hline
Wasserstein flows& No & No 
& $ \sum_{i=1}^{M/2} \int_{\CX} w_{2i} \phi(w_{2i - 1}x) \, dx$; (Measures) \\
\hline
Kernel SVM& Yes & No 
& $  w_0 + \sum_{i=1}^{N} w_i y_i K({\cdot, x_i}); M = N + 1$ \\ 
\hline
Finite Element 
& Yes & No 
& $  \sum_{i=1}^M w_{ij} \, \phi_i(\cdot) $; (Matrices with basis functions) \\
\hline
Feed-forward NN 
& No & No 
& $ \sum_{i=1}^W w^D_i \phi_i^D\!\left(\sum_{j=1}^W w_{ij}^{D-1} \phi^{D-1}_j(\cdot)\right)$; $M \sim DW^2$ \\
\hline
\end{tabular}

\end{minipage}

\caption{(\textbf{a}) Schema for methods, \hfill (\textbf{b}) Some typical parametrised architectures and the models produced by them.}
\label{fig:method_and_summary}
\end{figure}

\subsection{Architectures and their basic attributes}\label{arch_prec}

\begin{definition}[Architectures]\label{arch_defn}
We call $\CN:\R^M \to G, M \in \mathbb{N} \cup \{\infty\}$ an $M$ parameter model producing architecture (or simply architecture) if $\CN \in \cC^2_{pw}(\R^M, G)$\footref{notate}, with $\weightvector  \equiv [w_1, ..., w_M]  \in \R^M$ as its parameters. 
\end{definition}

\begin{remark}[Architecture features]\label{para_model_not} 
We denote the model set of $\CN$ as $G_M \eqdef \{\cNw: \weightvector \in \R^M\}$, the model subset near $\Phi$ as $\CBep^M \eqdef G_M \cap \CBep$, the set of optimal parameters as $\{\weightvector_M^{(*)}: \CL[\CN(\weightvector_M^{(*)})] \leq \CL[g],  \forall g \in G_M \}$, and:
\begin{equation}\label{para_model}
    \begin{aligned}
    & \CNw \eqdef \rmD \CN[\weightvector], \qquad\qquad\quad \, \theta(\weightvector) \eqdef \Ndwdag\Ndw, \qquad\qquad\qquad\qquad\qquad \CNww \eqdef \rmD^2 \CN[\weightvector] \\
    &\Thetaw \eqdef \Ndw \Ndwdag, \,\, \qquad \qquad \mu(\weightvector) \eqdef \inf \Big[\sigma \left( \thetaw \right)\setminus \{0\}\Big], \, \qquad\qquad\text{ where } \sigma(A) \eqdef \text{Spectrum of an operator } A \\
    \end{aligned}    
\end{equation}
We note again that in Defn. \ref{arch_defn}, for $M = \infty, \R^M = l^2$ or $\Rinf$. $\theta, \Theta$ are spectrally similar objects that encode how $G_M$ is \textit{immersed}\footref{immersions} and behaves in $G$ (Sec. \ref{arch_features}). We will comment on these aspects in Section \ref{arch_features}.
\end{remark}
\begin{remark}
Note the definition of an architecture is quite simple - this is what allows us to pack a large universe of methods within a single framework (Fig. \ref{fig:method_and_summary}\textbf{b}). Fortunately, most $\CN$ are piecewise regular over compact subsets of $\R^M$. Further, pointwise pathologies (as with ReLU) are irrelevant as we will ultimately be resolving more critical singularities than these and the class of $\cC^2_{pw}$ architectures is dense in the space of $\CN$ possible on a computer. 
\end{remark}

\begin{remark}
A property is “typical'' or “generic'' or holds “almost everywhere'' on an infinite dimensional set $X$, if it holds on a co-meagre subset of $X$ (Defn. \ref{Co_meagre}). In contrast, for finite dimensions, these notions may be supplied by the subset being co-meagre and/or having full Lebesgue measure, as appropriate.
\end{remark}

\begin{figure}
    \centering

    \begin{subfigure}[b]{0.479\textwidth}
        \centering
        \begin{tikzpicture}[
            neuron/.style={circle, draw, minimum size=18pt},
            weight/.style={midway, fill=white, inner sep=1pt},
            >=stealth
        ]

        \node[neuron] (x) at (-0.5,0) {$$};

        \node[neuron] (h1) at (3,2.15) {$\sin(\cdot)$};
        \node[neuron] (h2) at (3,-2.15) {$\sin(\cdot)$};
        \node (h3) at (2.95,0) {$\cNw(\cdot) = w_0 + \sum\limits_{i=1}^{2} w_{2i}\sin(w_{2i - 1} (\cdot))$};

        \node[neuron] (y) at (6.5,0) {$$};

        \node (bias) at (6.5,4.5) {$$};

        \draw[->] (x) -- (h1) node[weight] {$w_3$};
        \draw[->] (x) -- (h2) node[weight] {$w_1$};

        \draw[->] (h1) -- (y) node[weight] {$w_4$};
        \draw[->] (h2) -- (y) node[weight] {$w_2$};

        \draw[->] (bias) -- (y) node[weight] {$w_0$};
        \end{tikzpicture}
    \end{subfigure}
    \hfill
    \begin{subfigure}[b]{0.479\textwidth}
        \centering
        \begin{overpic}[width=\textwidth]{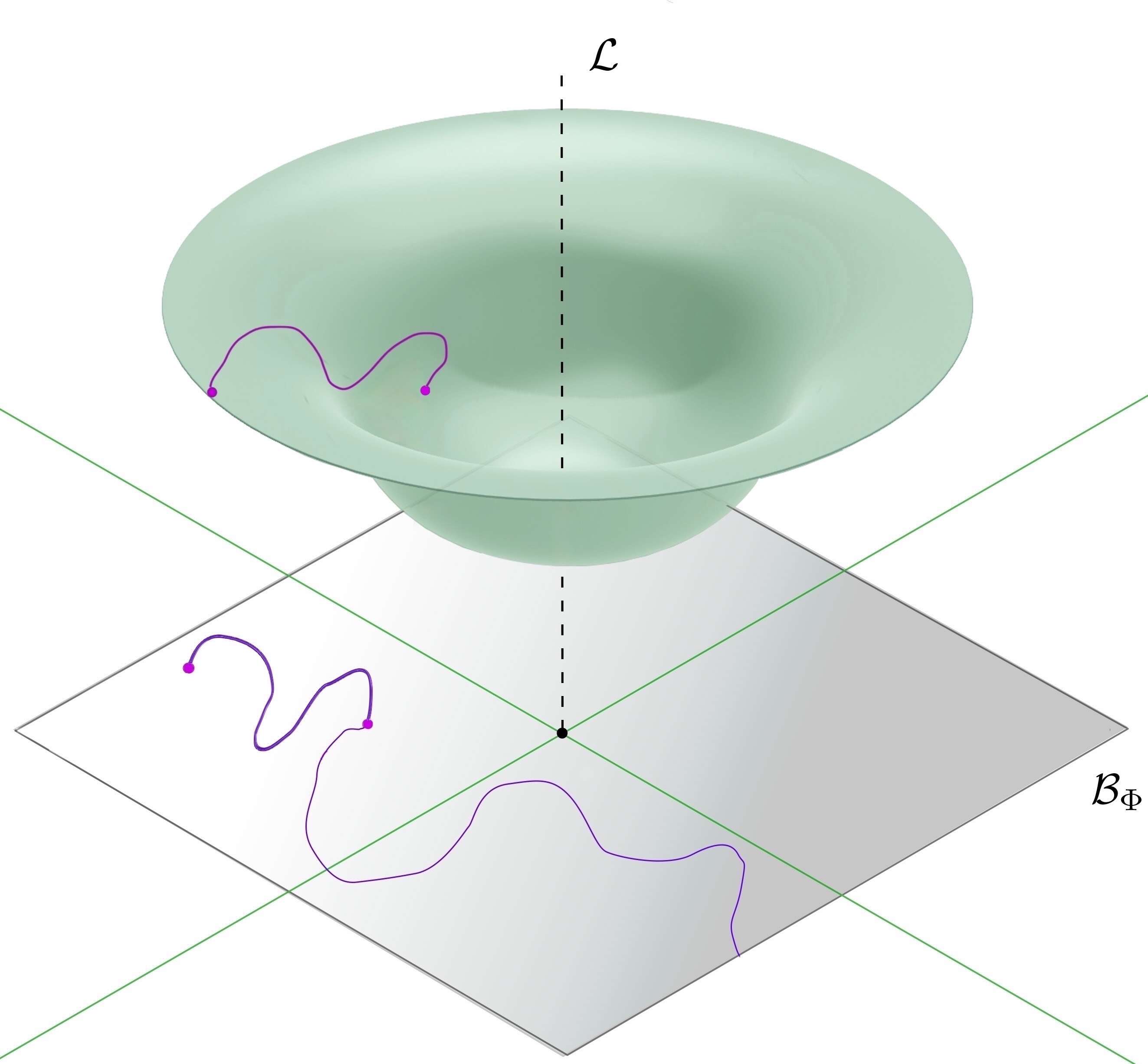}
            \put(53,29){\makebox[0pt][s]{{$_{\dot{\CN}_1(t) \, = \, -\Theta_1(t) \nabla\CL[\CN_1(t)]}$}}}
            \put(22,38){\makebox[0pt][s]{{$_{\CN_1(t)}$}}}
            \put(63,33){\makebox[0pt][s]{{$_{\CN_1(0) \, = \, g(0)}$}}}
            \put(39.5,18.3){\makebox[0pt][s]{{$_{\leftarrow \, \CBep^1 \, \to}$}}}
        \end{overpic}
    \end{subfigure}
    \caption{(\textbf{a}) NN with $a = 2, M = 5$, and sine activations, \hfill  (\textbf{b}) Parametric gradient flow dynamics}
    \label{NN_and_param_Gd}
\end{figure}

Fourier series fits, deep Neural Networks (NNs), Large Language models, ensembles of grid cell neurons, etc, all fit this characterization. Our stock class of examples is comprised of Neural Networks (Fig. \ref{NN_and_param_Gd}\textbf{a}), real-analytic and non-linear w.r.t. $\weightvector$, and has a dense model set $\{\cNw: \weightvector \in \R^M\} \subset L^2([-\pi,\pi])$ if $M = \infty$. Ex. \ref{stock2} is a generalized Fourier series estimator, with tunable co-efficients \textbf{and} frequencies. 

\begin{example}[Fig. \ref{NN_and_param_Gd}\textbf{a}]\label{stock2}
Let $\weightvector \in \R^{M}, G = W^{2,2}([-\pi,\pi])$. Our example architecture $\CN$ is given by
\[\CN: \R^M \to G, \quad \cNw(\cdot) = w_0 + \sum_{i = 1}^{a} w_{2i}\sin \Big(w_{2i-1}(\cdot) \Big), \quad M = 2a+1\]
\end{example}

Having formally setup architectures, we can now parametrize $G$ using Defn. \ref{arch_defn} and get feasible gradient flows.

\subsection{Parametric optimisation and its connection to gradient flows}\label{Optimization_Flows}

We now need a notion that translates the nominal loss and problem into parametrised forms, using the parameters and generated models of an architecture $\CN$. Further, we would want to enforce that optimisation under these parametrisations leads to better models. Fundamental notions from Lyapunov and gradient theory are needed here. 

\begin{definition}[Parametric Loss and Gradient Flows]\label{generic_para}
Let $\CN$ be an architecture. We call $\cL \in \cC^2(\R^M \times G_M)$ a parametric loss function for $\CL$, if $\CL$ is a strict Lyapunov function\footref{lyapunov} for the following dynamical system:
\[\frac{d}{dt}\wt = -\nabla \cL[\cNw, \weightvector], \qquad \qquad \qquad \weightvector \in \R^M, \cNw \in G_M
\]
\end{definition}

\begin{remark}[Update/Optimisation Rule]\label{loss_para}
Unless stated otherwise, for a given $\CL$, we shall define the parametric loss for $\CN$ by $\cLw \eqdef \CL[\CN(\weightvector)]$ and associated update/optimisation rule by the following gradient flow:
\begin{subequations}\label{param_GD}
\begin{alignat}{2}
\frac{d }{dt}\wt = -\nabla \cL[\wt] = -\Ndwdag \nabla\CL[\CNt], \qquad\qquad\qquad \wt\in\R^M \label{param_GD1}
\\
\implies \frac{d }{dt}{\CNt} = -\Ndw\nabla \cL[\wt] = -\Thetat\nabla\CL[\CNt], \qquad\quad\textbf{ } \CNt \in G_M \label{GMdynamics}
\end{alignat}
\end{subequations}
\end{remark}
\noindent where we let $\CN(\wt) \equiv \CNt, \Theta(\wt) \equiv \Theta(t)$, etc. Note that \textbf{all} ODEs with a strict Lyapunov function $f$ are effectively a gradient flow w.r.t. $f$ under the correct coordinate chart \cite[Thm. 1]{Barta12}. Thus, if our optimisation method has $\CL$ (and therefore $\CL \circ \CN$) as a Lyapunov function (a natural expectation from any differentiable optimisation schema), the choices made in Remark \ref{loss_para} are extremely well justified vis a vis analysis. 

We are now able to formalise the notion of having a \textit{good initial model} or more precisely, being well-initialized:

\begin{hypothesis}[Well initialised]\label{well_initialize}
Let $\,\, \frac{d}{dt}\wt = \CF[\cNw, \weightvector, t]  \,\,$ represent some evolution equation over the parameters of $\CN$. We say $\CN$ is well-initialized under the evolution equation with parameters $\wzero$, if there exists $t \in \R^+$ s.t. $\CNt \in \CBep$. $\CN(0) \equiv \CN(\wzero)$ is called a well-initialized model.  
\end{hypothesis}

\begin{remark}\label{GD_analysis_is_sufficient}
    Note the optimization flows can be imposed using more complex rules than Eq. \ref{param_GD}. These choices are where data often plays a central role. As Hypothesis \ref{well_initialize} is generic, our discussions ahead will strictly lower bound the expected performance of the “optimal'' choices. 
\end{remark}

Well-initialisation may seem like an unreasonably strong hypothesis but methods like stochastic gradient descent (SGD) imply we can often obtain such models if $\CBMep$ is non-empty \cite[Thm. 1]{sgd_global_convg_geman_86} (see App. \ref{SGD_uses_lin_vs_nlin_N}). However, while excellent at obtaining initial guesses, convergence rates of these methods are usually too slow post the \textit{initial search} phase. Further, “good guesses'' are often validated “heuristically'' - we do not prescribe formal rules for validating models obtained this way. Regardless, the convergence phase needs better understanding and algorithms, specially when well-initialisation is a practical reality with coarse-grained guesses \cite{ good_guess_k_Means_13, good_guess_then_refine_Ilya_13,warm_start_Yildrim_02}.

However, before we can establish the parametric equivalent of Lemma \ref{simple_GD_G}, we have to intuit why nonlinear $\CN$ complicate analysis and why they are still used when linear $\CN$ are much simpler and nominally sufficient.
\vspace{\baselineskip}

\subsection{Linear vs Nonlinear Architectures}\label{nlin_A}

\begin{figure}[h]
\centering
\begin{subfigure}{0.48\textwidth}
    \begin{overpic}[width=\textwidth]{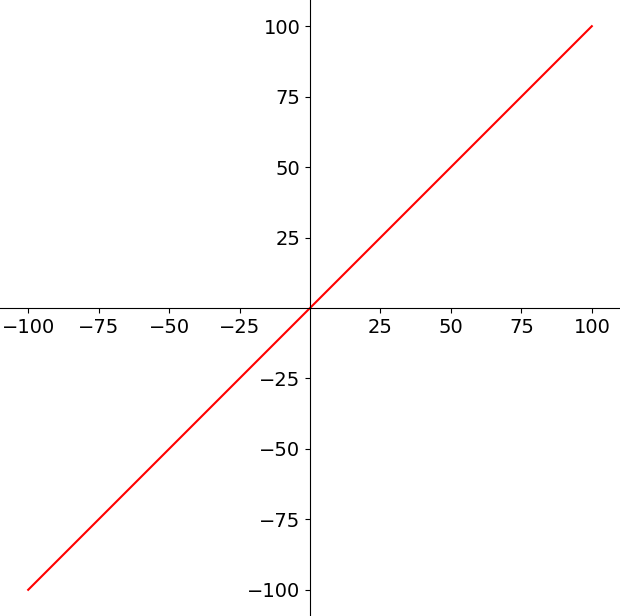}
        \put(10,90){$G = \R^2$}  
        \put(60,10){$\CN_l (\weightvector) = (\weightvector, \weightvector)$}
    \end{overpic}
    \label{fig:Lin_vs_NLin_first}
\end{subfigure}
\hfill
\begin{subfigure}{0.49\textwidth}
    \begin{overpic}[width=\textwidth]{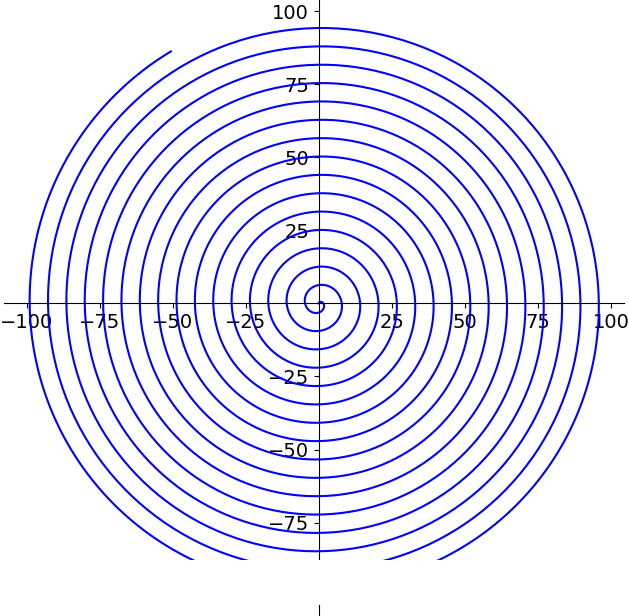}
        \put(10,90){$G = \R^2$} 
        \put(20,4.5){$\CN_n(\weightvector) = ({\weightvector}\sin(\weightvector), {\weightvector}\cos(\weightvector))$}
    \end{overpic}
    \label{fig:Lin_vs_NLin_second}
\end{subfigure}        
\caption{A linear 1-parameter architecture vs a nonlinear 1-parameter spiral architecture.}
\label{fig:Lin_vs_NLin}
\end{figure}

For all architectures, Picard-Lindelof theorems guarantee the solutions $\wt,\CNt$ exist for Eq. \ref{param_GD} for all $t \in \R^+$. Further, if $\FF$ is well-behaved, $\CN$ is linear w.r.t. $\weightvector$, and $M < \infty$, then $\cL$ is a well-behaved map and Eq. \ref{param_GD} comprises a well-understood paradigm for studying optimization \cite{Convex_Opt_04}\cite[Ch. 4]{manyfold_learning}. Thus, linear $\CN$ are always sufficient if $d(\Phi, G_M)$ is \textit{small}. However, linear $\CN$ access $\leq M$ dimensional subspaces and cannot efficiently cover a truly $> M$ dimensional $\CBep$ (see Fig. \ref{fig:Lin_vs_NLin}). Unless $\Phi$ is known in advance, there is no way to \textit{a priori} force $d(\Phi, G_M)$ to be small.

In contrast, even for strongly well-behaved $\FF$, if $\CN$ is nonlinear and not injective (a usual occurrence for ML architectures), Eq. \ref{param_GD1} and \ref{GMdynamics} are on unequal footing and it is unclear when $\cL$ is well-behaved: the optimisation regimes are of relatively limited understanding. However, using $M$ parameters nonlinearly, we can efficiently \textit{fold} $G_M$ in $G$ so $\CBMep$ is \textit{dense enough} ($d(\Phi, \CBMep)$ is small), even if $M < \dim(G)$ \cite[Sec. 3]{Pinkus_bounded_depth_widht_99}. Unfortunately, the better coverage comes with twin costs: $\cL$ can have too many spurious critical points and/or pathological regions of ill-behavior for convergence under Eq. \ref{param_GD}, even if $\CN$ is well-initialised and the nominal problem is well-behaved.

To exemplify, assume $G= \R^2, \CBep = [-100,100]^2$, and $\CL[g] = \braket{g - \Phi, g - \Phi}$ for some
$\Phi \in \CBep$. Consider some fixed linear $\CN_l:\R \to \R^2$ and a nonlinear $\CN_n: \R \to \R^2$ given by $\CN_n(\weightvector) = ({\weightvector}\sin(\weightvector), {\weightvector}\cos(\weightvector))$ as in Fig. \ref{fig:Lin_vs_NLin}. For any $\Phi$, the associated optimal $\CN_l(\wst)$ will always be found using GD on $\cL_l = \CL[\CN_l]$, but $\CN_l$ is very limited in its coverage of $\CBep$. In contrast, the nonlinear $\CN_n$ covers $\CBep$ significantly better, but unless $\wzero$ is near enough the global minimum of $\cL_n$, even if $\CN_n(0) \in \CBep$, $\CNt$ will still get stuck in a sub-optimal minimum (usually, $\dim(G) = \infty$ - Fig. \ref{fig:Lin_vs_NLin} and this argument are only motivating how nonlinearity can be advantageous and disadvantageous). 

The density of $\CBMep$ in $\CBep$ controls the best we can do and the geometry of $G_M$ in $G$ guides where $\CNt$ will go. The need to precisely understand such facets of parametric optimisation leads us to the last component of learning.

\section{Setup, convergence, and model set analysis (Mathematics)}

\begin{figure}[ht!]
\centering

\begin{minipage}{0.3\textwidth}
\centering
\begin{tikzpicture}[
    every node/.style={
        draw, 
        circle, 
        minimum size=1cm,
        inner sep=2pt,
        align=center
    },
    >=Stealth,
    thick
]

\node (A) at (135:2) {Properties\\ of $\FF, G, \Phi,$\\ and $\CL$};
\node (B) at (45:2) {Properties\\ of $\CN, G_M,$\\ and $\cL$};
\node (C) at (225:2) {Dynamics\\\& limits\\ in $G_M$};
\node (D) at (315:2) {Pruning,\\extending\\$\CN$};

\draw[<->, bend right=25] (B) to (A);
\draw[<->, bend right=25] (A) to (C);
\draw[<->, bend right=25] (C) to (D);
\draw[<->, bend right=25] (D) to (B);

\node[draw=none] at (0,0) {\textbf{Analysis}};

\draw[<-] (C) -- +(1.15cm,1.15cm);
\draw[<-] (B) -- +(-1.15cm,-1.15cm);
\draw[<-] (A) -- +(1.15cm,-1.15cm);
\draw[<-] (D) -- +(-1.15cm,1.15cm);

\end{tikzpicture}
\end{minipage}
\hfill
\begin{minipage}{0.65\textwidth}
\centering
\renewcommand{\arraystretch}{1.15}
\small
\begin{tabular}{|l|l|l|}
\hline
\textbf{Analysis} 
& \textbf{Focus} 
& \textbf{Key concepts and subfields} \\
\hline
Problem Setup [\ref{well-behaved}] & $\FF, G, \CL$ & Banach Manifolds, Gradient flows\\
\hline
Convergence in $G$ [\ref{well-behaved}] & $\CL, G, \Phi$ & Lojasiewicz Inequality (LI) \\
\hline
Method Setup [\ref{arch_defn}] & $\CN, \R^M$ & Lyapunov Theory, Gradient Flows \\
\hline
Optimisation flow [\ref{Optimization_Flows}] & $\nabla\CL, \Theta$ & Pullback metrics, NTKs \\
\hline
Well-initialisation [\ref{Optimization_Flows}] & $\CN, \cL$ & Stochastic Gradient flows \\
\hline
Density of $G_M$ \cite{GKP20, UAP_Wasserstein_space_Lu_20} & $G, G_M$ & Universal Approximation Theorems \\
\hline
Efficiency of $G_M$ [\ref{arch_features}] & $\CN, G_M$ & Immersion Theory, Transversality \\
\hline
Convergence in $G_M$ [\ref{FNtogether}] & $\CN, \cL$ & Desingularization, Fredholm theory \\
\hline
Expanding $\CN$ [\ref{expand_Arch}] & $\weightvector, \phi, \CN$ & Direct Limits, Extension theorems  \\
\hline
Pruning $\CN$ [\ref{prune}] \cite{reddog22} & $\weightvector, \phi, \CN$ & Renormalisation Group, Heuristics \\
\hline
\end{tabular}
\end{minipage}
\caption{(\textbf{a}) Analysis is interconnected, $\hfill$ (\textbf{b}) Types of analysis and concepts needed for different facets of learning.\hfill}
\label{fig:anaysis_breakdown}
\end{figure}

Analysis is expectedly prevalent in all kinds of learning - and thus, learning about learning. We have already indulged in some forms to setup our problems and methods, speak of convergence to solutions in nominal terms, etc. Many of our choices have been forward looking: we picked the Banach-Riemannian setting for $G$ because we wanted gradients, we used the Lojasiewicz Inequality because we wanted a weak condition for convergence guarantees, etc. If we may call this inferred or implicit analysis, it is time to consider the explicit and deductive forms - the analysis that follows after we have made some or all of our choices, but before the actual parametric optimisation.

Within the goal of establishing unifying results for learning, we may need to answer the questions concerning: the existence of an appropriate $\CL$ and $\Phi$ (to establish the problem is solvable), the density of $G_M$ in $G$ (to establish $\CN$ is powerful enough, at least in principle), convergence of parametric optimisation, etc (Fig. \ref{fig:anaysis_breakdown}\textbf{(a)}). There are several benefits to such anticipatory analysis. For example, convergence/rates established in \cite{feydy_20} were key to building the shape recognition tools in \cite{dog23}. Similarly, in \cite{cohen23, dog23_3}, establishing prior results on $\FF, \CN$ directly offered strategies for building differential equation solvers that significantly out-performed the baseline methods they were based on.

We do not establish novel $G_M$ density results, since Universal Approximation Theorems on various choices of $G$ already cover that. For example, a sequence of works has established that many common ML activations $\phi \in \cC^2(\R)$, in combination with a countable number of parameters, generate a set of Neural Network models that can be dense in the space of $\cC^k$ functions \cite{Cyb89, Bouned_Relu_SHEN_22}, in Sobolev spaces \cite{GKP20}, in Wasserstein \cite{UAP_Wasserstein_space_Lu_20} and Bochner \cite{UAP_Bochner_space_Schmocker_23} spaces, and in the space of continuous nonlinear operators between Banach spaces \cite{UAP_Mishra_21, Deeponet_error}: $G_M$ density is an essentially settled topic.

However, we can aim to augment these by invoking additional manifold like structures on $G_M$ to better study the dynamics possible on them. As it happens, the properties of $\theta, \Theta$ are critical in pretty much every aspect of model set analysis. The following results thus become key to understanding learning as a unified discipline:
\begin{itemize}
    \item Efficiency of $G_M$ in $G$: ensuring $G_M$ is as non-degenerate as possible to maximise the use of each parameter.
    \item Convergence in $\R^M$ and $G_M$: how models may inherit convergence properties from a well-behaved problem.
    \item Expanding architectures without performance losses: refining a good initial model by adding more parameters.
\end{itemize}
As we noted before, $\theta,\Theta$ are key to the geometry of $G_M$ in $G$ and thus to answering such questions. Let us see how:

\subsection{When are model sets efficient and how do architecture features decide that?}\label{arch_features}

Let us explicitly write out $\theta(\weightvector), \Thetaw$ for a better perspective on how they may be viewed computationally:
\begin{equation}\label{theta_defn}
\Ndw = \begin{pmatrix}
\partial_{w_1} \cNw \\
: \\
\partial_{w_M} \cNw \\
\end{pmatrix}, \qquad \thetaw = \Ndwdag\Ndw \implies \theta_{ij}(\weightvector) = \Bigg\langle{\frac{\partial \cNw}{\partial w_i}, \frac{\partial \cNw}{\partial w_j}}\Bigg\rangle 
\end{equation}
In turn, we may represent $\Thetaw$ in coordinate free terms by its action on $g \in G$ as:
\begin{equation}\label{Theta_defn}
\Thetaw g = \sum_{i=1}^M \Bigg \langle \frac{\partial \cNw}{\partial w_k} , g \Bigg \rangle \frac{\partial \cNw}{\partial w_k}    
\end{equation}

Architectures $\CN \colon \R^{M} \to G_M \subset G$ are differentiable maps into a Banach manifold by construction. Drawing on notions and abusing notation from differential geometry, we shall sometimes view the image of  $ \Ndw \colon \R^M \to G$ as the tangent space of $G_M$ at $\CN(\weightvector)$. However, we emphasise that although this is indeed a linear subspace of $G$, identifying it as a tangent space is strictly a formal calculation. 
The adjoint of this map $\Ndwdag \colon G \to \R^M$ can in a similar fashion be viewed as a ``projection'' onto this tangent space. We can then pull back using the loss from $G$ to $\R^M$ to generate the parameter flow (Eq. \ref{param_GD1}) and view $\Theta$ as projecting $\nabla \CL$ onto $G_M$ for the model flows (Eq. \ref{GMdynamics}). 

$\theta,\Theta, \mu$, etc are significant because they are the computable controls on the geometry of $G_M$ in $G$. Our treatment and interpretation of $\theta, \Theta$ is novel in some key ways, but the uses will be similar to existing ones in the literature: indeed, the ``NTK'' \cite{jacot_ntk_18}, Koopman operator \cite{dog20, dogred20}, ``pull-back metric'', etc, refer to related objects that differ only due to the context they are invoked in \cite{dog25_6}, but perform similar roles. In particular, our approach considers $\weightvector$ as the ``input'' w.r.t. which the NTK is being considered, not the data. Data is considered an extrinsic driver of the dynamics, separate from the intrinsic geometry $G_M$ takes in $G$ once $\CN$ is chosen and the gradient flows Eq. \ref{param_GD} imposes. For us, this geometry is the key object of study and $\Theta$, as defined in Defn. \ref{arch_defn}, our ``NTK''.

A common and critical assumption in conventional NTK theory is the invertibility of $\Theta$ \cite{jacot_ntk_18, cosine_ntks_full_rank_prob} (within the context it is invoked in). A related point of interest is in tracking the spectrum of $\Theta$ since that allows control on optimisation rates. Thus, it is interesting to observe that $\theta, \Theta$ share their non-zero spectrum and that typical $\CN$ are immersions\footref{immersions} and their $G_M$ are immersed sub-manifolds in $G$ almost everywhere. That is, intuitively speaking, none of the $M$ local directions collapse near $\cNw \in G_M$ and $G_M$ possesses a \textit{mostly} manifold like structure locally:
\begin{lemma}[Typical architectures produce manifold-like model sets \cite{manyfold_learning}]\label{imm_sub_man_full_rank}
    Let $\CN$ be an architecture and $\theta, \Theta$ be defined as per Defn. \ref{para_model_not}. Then we have $\sigma \left( \thetaw \right)\setminus \{0\} = \sigma \left( \Thetaw \right)\setminus \{0\}$. Further, for a generic architecture $\CN$, for $\weightvector \in \R^M$ almost everywhere, $\theta$ is full rank and $G_M$ is locally an $M$-dimensional immersed submanifold\footref{immersions} at $\cNw$. 
\end{lemma}

\begin{remark}
Through Lemma \ref{imm_sub_man_full_rank}, we have that large classes of architectures allow $\thetaw$ to be full rank a.e. in $\R^M$ \cite{cosine_ntks_full_rank_prob}. Note that linear $\CN$ can easily produce $M$ dimensional $G_M \subset G$: the utility of a nonlinear $\CN$ without full-rank $\theta$ a.e. over the parameter space $\R^M$ is usually low and would make for an inefficient choice for parametrising $G$. 
\end{remark}

Thus, nonlinear optimization is as much a dynamics on an a.e. immersed submanifold $G_M$ in the target space $G$, as it is a dynamics in the parametric space $\R^M$. Unsurprisingly, $\theta,\Theta$ variants have been used to characterize $\CNt$ dynamics in specific cases, such as for models trained on finite datasets \cite{Loja_to_ML_2022} and over-parametrised $\CN$ \cite{jacot_ntk_18}, etc.

Lemma \ref{imm_sub_man_full_rank} also implies that $\Theta$ is an essentially $M$ dimensional operator and significantly simplifies the complexity of analysing the model dynamics implied by Eq. \ref{GMdynamics}. For example, its spectrum can be directly tracked from the spectrum of the $M \times M$ matrix $\thetaw$, which is always estimable using Eq. \ref{theta_defn} when a gradient step is computed.

A well-behaved $\CL$ guaranteed the existence of tractable dynamics over $G$: well-behaved $\theta, \Theta$ will allow us similar statements on $\R^M, G_M$. For example, when $\CN$ is linear, the difference between Eq. \ref{param_GD} and Eq. \ref{idealized_G_flow} is simply the constant $\Theta$ term in front of $\nabla \CL$: as long as we have control on $\Theta$, we have control over the optimization dynamics.

Analytic $\CN$ will be a useful example of such a class: for example, $\thetaw$ is full rank \textit{a.e.} in $\R^{M}$ for Ex. \ref{stock2}. If we were solving Ex. \ref{stock} (a strongly well-behaved problem) with Ex. \ref{stock2}, this imbues $(\wt, \CNt)$ with determinable convergence properties w.r.t. Eq. \ref{param_GD}  (Thm. \ref{vtk_central}). However, to formally support these claims, we must answer:
\begin{itemize}
    \item When is $\cL$ well-behaved around its critical points? (a result like Lemma \ref{simple_GD_G} on $\cL$)
    \item If $\CN$ is nonlinear and $\cL$ well-behaved, when is convergence locally and/or globally optimal w.r.t. $G_M$?
    \item If $\Phi \notin G_M$ for any finite $M$, how do we find models limiting to $\Phi$? 
\end{itemize}

\subsection{Convergence in parameter spaces and model sets}\label{FNtogether}

We recover some geometric structure and convergence guarantees by showing that while $\CNt \in \CBep$ may not guarantee we will eventually reach $\Phi$, it can lead to some $\CNwst \in \CBMep$ at determinable rates. We begin by characterizing the critical points of the flow $\wt$ within its omega-limit set\footref{omega_limit_sets} $\omega(\weightvector)$ (limit points of future iterations):
\begin{lemma}\label{equilibria}
    Assume $\FF$ is well-behaved, $\CN$ is well-initialised, and $\wst \in \omega(\weightvector)$ is a critical point of $\cL$ w.r.t. the system prescribed by Eq. \ref{param_GD}. Then, one or more of the following holds:
    \[
    \textbf{(i) } \nabla \CL[\CN(\wst)] = 0, \qquad\qquad \textbf{(ii) }
    \Theta(\wst) = 0, \qquad\qquad \textbf{(iii) } \Theta(\wst) \neq 0 \neq  \nabla \CL[\CN(\wst)] \in \ker(\Theta(\wst))    \]
\end{lemma}

Lemma \ref{equilibria} tells us that parametric models stop improving only if: \textbf{(i)} they are already at $\Phi$, \textbf{(ii)} the model set $G_M$ becomes singular (the manifold gets pinched to a point like structure), or \textbf{(iii)} the desired direction of the flow is orthogonal to the model set. Lemma \ref{imm_sub_man_full_rank} already established that condition \textbf{(ii)} is rare. Further, $\Theta$ and $\nabla\CL$ are independent objects - there is no reason why $\Thetaw$ should be singular exactly where $\nabla \CL$ guides the flow.

In contrast, since $G$ is often infinite dimensional and $G_M$ is usually finite dimensional, condition \textbf{(iii)}, the  orthogonality of $\Theta$ and $\nabla \CL$ is an ever-present concern for optimisation. Consider Ex. \ref{stock2} again: $\Theta$ is never 0, but under a mean squared $\CL$, there are too many equilibria simply because $\Thetaw$ and $\nabla\CL[\cNw]$ would be orthogonal at those parameters. We have a twin issue: we need to understand the optimisation flows and their stopping points in these terms and we need to figure out a way to restart the flow once we encounter such equilibria. 

We can resolve both by specifying when $\cL$ is locally well-behaved around its equilibria. While the proof needs more setup, the geometry and dynamics of parametrised models is fully expressed in terms of the introduced objects.

\begin{theorem}[Ch.4, \cite{manyfold_learning}]\label{vtk_central}
Assume $\FF$ is well-behaved w.r.t. $\CL$ and $\CN$ is well-initialised under Eq. \ref{param_GD}. Then, $\cL = \CL \circ \CN$ is locally well-behaved at a critical point $\wst \in \omega(\weightvector)$, if any of the following hold for all $\weightvector \in \CBdw$:
    \begin{enumerate}
        \item $\exists \alpha^* \in [\frac{1}{2}, 1), C^* > 0$, s.t. $|\cL[\weightvector] - \cL[\wst]|^{\alpha^*} \leq C^*\|\DcLw\|$.
        
        \item $M \in \N$ and $\cL$ is definable in an o-minimal structure\footref{o_minimal_structure}. 

        \item $\mu(\weightvector) \geq c > 0$. In this case, $\alpha^* = \alpha$. 
        
        \item For some $r > 0$, $\mu(\weightvector)$ satisfies the distance LI: $|\mu(\weightvector)| \geq |\weightvector - \wst|^r$. In this case, $\alpha^* = \frac{\alpha + r}{1 + r}$.
    \end{enumerate}
\noindent Further, we have $\wt \xrightarrow[t \to \infty]{} \wst \in {R^M}$ with convergence rates:
\begin{equation}\label{LI_convg_rates}
\|\wt - \wst\| = \begin{cases}
                O(e^{-C^*t}), \text{ if } \alpha^* = \frac{1}{2}\\
                O(t^{-\frac{1 - \alpha^*}{2\alpha^* - 1}}), \text{ if } \frac{1}{2} < \alpha^* < 1\\
                \end{cases} 
\end{equation}
Finally, $\cL$ is well-behaved if it satisfies one of these conditions at its global minima $\weightvector_M^{(*)}$.
\end{theorem}
\begin{remark}
$(1)$ imparts well-behavior on $\cL$ directly by Hypothesis \ref{well-behaved}. $(2)$ generically implies $(1)$ by the celebrated work of Lojasiewicz \cite{Loja_og_63, Loja_og_65, Loja_og_83} on analytic functions that was then generalized to functions definable in an o-minimal structure by Kurdyka \cite{kurdyka98} (see \ref{o_minimal_structure} for precise definitions - intuitively, the graph set of $\cL$ is geometrically “simple'' and does not contain pathologies like a Cantor set, strong singularities, etc). $(3)$ implies $(1)$ because $\CL$ already satisfies \ref{Lojasiewiczineqaulity} over $\CBep$ and $\CN$ is effectively a diffeomorphism (modulo some subspaces that we can quotient out) and contains sufficient structure to pass \ref{Lojasiewiczineqaulity} on to $\cL$. $(4)$ implies $(1)$ via some fundamental de-singularisation arguments.

Conditions $(1, 2)$ are not the focus of our work from a theoretical perspective though we may use them for applications. There is a large literature that exists on the intersections of studying “tame'' geometries/topologies and their connections to \ref{Lojasiewiczineqaulity} type relations that the reader may consult \cite{bier_milman_semianalytic_88, Loja_Chill_inherits_LI_03, kurdyka98}. The propagation of \ref{Lojasiewiczineqaulity} type relations on functionals $\CL$ from the Hilbert space $G$ to restrictions on subsets/submanifolds $G_M$ within them and their applications comprise a more novel line of enquiry, with some precedent in \cite{res_singular_feehan_19, Loja_on_submanifolds_Rupp_20}. We are unaware of prior work studying the propagation of \ref{Lojasiewiczineqaulity} to parametrised sub-manifolds to understand what implications pull-backs from the target space to the parameter space have for parametric and model optimization dynamics. 
\end{remark}
\vspace{\baselineskip}

We list some immediate corollaries of Thm. \ref{vtk_central} that are relevant for commonly used $\CN$:
\begin{corollary}\label{example-Li}
Assume $\FF$ is well-behaved and $\CN$ is well-initialised. $\cL$ is locally well-behaved at a critical point $\wst \in \omega(\weightvector)$, if any of the following conditions hold:
\begin{itemize}
    \item $\cL$ is analytic on $\CBdw$ and $M \in \mathbb{N}$.

    \item $\CN$ is affine on $\CBdw$ and $M \in \mathbb{N}$. In this case, $\alpha^* = \alpha,  C^* = \mu(0)C$. 

    \item $\CN$ is an immersion on $\CBdw$ and $M \in \N$. In this case $\alpha^* = \alpha$. 

    \item $\CN$ is analytic on $\CBdw$ and $M \in \mathbb{N}$. In this case, $1 > \alpha^* \geq \alpha$.
\end{itemize}
\end{corollary}

Thm. \ref{vtk_central} and Corollary \ref{example-Li} break down the \textit{computational tractability} of a \textit{nominally solvable} problem into a question on the parametrised models we use (for example, the analytic condition is satisfied in Ex. \ref{stock2}). The four classes, with their resulting corollaries are essentially a question of being able to force \ref{Lojasiewiczineqaulity} directly or showing that $G_M$ is geometrically \textit{tame} enough to conserve the structure of the well-behaved $\CL$.

Note that $\cL$ being well-behaved is not as powerful for $\cL$ as it was for $\CL$: $\CNzero \in \CBep$ would converge to $\Phi$ if the dynamics were in $G$, but since Eq. \ref{param_GD} operates over $\R^M$ and $G_M$, there are usually too many potential equilibria $\wst$ around which $\cL$ is locally well-behaved: we have no guarantee of reaching the global minima $\weightvector_M^{(*)}$. Qualitatively, Theorem \ref{vtk_central} only establishes refinements at \textbf{determinable} rates to a ``local'' solution $\CN(\weightvector^*) \in \CBMep$. 

However, we still establish qualitative behavior and useful special cases. For example, $\Thetat$ in Eq. \ref{param_GD} acts as a constraint operator at all $t$, forcing $\CNt$ to stay in $G_M$, even if $\nabla\CL$ attempts to push it out to generic elements in $G$. Thus, Eq. \ref{param_GD} with a constant $\Thetat$ (effectively affine $\CN$, such as the conventional NTK regime) is equivalent to the nominal flow over $G$ (Lemma \ref{simple_GD_G}), with the globally optimal model $\CN(\weightvector_M^{(*)})$ playing the role in $G_M$ that $\Phi$ plays in $G$, if $\CL$ is at least strictly well-behaved. This is the essential idea behind why linear $\CN$ can be so powerful.
 
For nonlinear $\CN$, the usual non-trivial task is understanding how $G_M$ is ``immersed'' in $G$. Well-behaved $\CL$ on $\CBep$ provide only part of the needed structure - we may have regained some convergence properties, but we don't yet know where we go. Indeed, even for strongly convex $\CL$, $d(\CNt, \Phi)$ can increase as $\cL(t)$ decreases.

Further, while infinite dimensional well-behaved problems are in principle solvable (with infinite dimensional gradients - classical NTK theory in essence), computationally, we are limited to finite dimensional methods without such guarantees. Nonlinear $\CN$ are especially saddled with an abundance of critical points and an approach to the ``best'' model possible within $G_M$ can't be guaranteed (Fig. \ref{fig:Lin_vs_NLin}).

Thus, even though unravelling how $G_M$ is ``immersed'' in $G$ provides greater control via Thm. \ref{vtk_central}, we still don't have a method for arbitrarily good $\Phi$ estimates. We remedy this situation by providing a construction that bypasses this in a computationally realizable way, which is the chief limitation of the classical NTK perspective \cite{Mariia_NTK_Collapse_23}. 
\vspace{\baselineskip}

\subsection{Architecture expansions and sparsification}\label{arch_exp_prun}

Historically, architectures were built to work with the least number of parameters possible, owing to the laborious human calculations involved and principles like Occam’s razor. Until recently, computational modelling methods also adhered to those philosophies since larger models cost more. In contrast, machine learning (ML) often leans heavily in the opposite direction, with the advent of
powerful computers (and GPUs) having spurred a philosophy of “more is better” \cite{Moreisbetter_simon24}. However, sparsification/pruning of architectures (modifying $\CN$ via removal of redundant parameters to save on computational costs) remains a highly active sub-area due to the high compute costs \cite{burkholz2021existence, ele21, mae21, tanaka2020pruning}. 

Thus, the principled expansion and pruning of architectures are important facets of parametrised learning. Poetically, they can be seen as two sides of a single motivation: to change $\CN$ \textbf{without} losing performance. We want the simplest model set geometry that can support the optimisation dynamics as efficiently and long as needed.

For example, often $\dim(G) = \infty$. We cannot usually ensure $G_M$ contains $\Phi$ for some $\CN$ with $M \in \N$: it follows that $\Phi$ can only be approximated arbitrarily well if we can augment $\CN$ by adding more parameters. One useful philosophy might be to “\textbf{learn, then expand}'', instead of simply optimising over the largest feasible $\CN$ (expand, then learn). But how do we add parameters without losing already obtained performance and when do we do so? 

A natural idea would be to  expand $G_M$ and access new gradients at local equilibria $\CNwst$, without losing the model $\CNwst$ during the expansion. Ideally, expansions would fit our geometric viewpoint in a manner that recovers the “over-parametrised limit'' results too \cite{ntk_imp_19, jacot_ntk_18}, which are mathematically simpler but computationally infeasible.

Similarly, we would want $\CN$ pruned in a manner that reduces $M$ as much as possible \textbf{without} impeding the optimisation dynamics. A natural choice would be the non-equilibria points $\cNw$ where $G_M$ becomes degenerate \textbf{without} pausing the flow. Ideally, this would match the established theory on pruning \cite{fra20, reddog22}.

\newpage
Let us set some basic notation that is handy in considering these matters and begin with expansions:

\begin{remark}[Parametric projections and extensions]\label{para_projection_extensions} 
Let ${M_j} > {M_i}$. $\{w_1, ..., w_{M_i}\}$ is taken as the canonical $\proj_{\R^{M_i}}(\weightvector_j)$ of $\weightvector_j \equiv \{w_1, ..., w_{M_j}\} \in \R^{M_j}$. The canonical extension of $\weightvector_i \equiv \{w_1, ..., w_{M_i}\} \in \R^{M_i}$ to $\R^{M_j}$ is $ \{w_1, ..., w_{M_i}, 0, ..., 0\} \in \R^{M_j}$ and often denoted as just $\weightvector_i$. $\wi^*$ and $\wi^{(*)}$ will represent equilibria and global minima.
\end{remark}

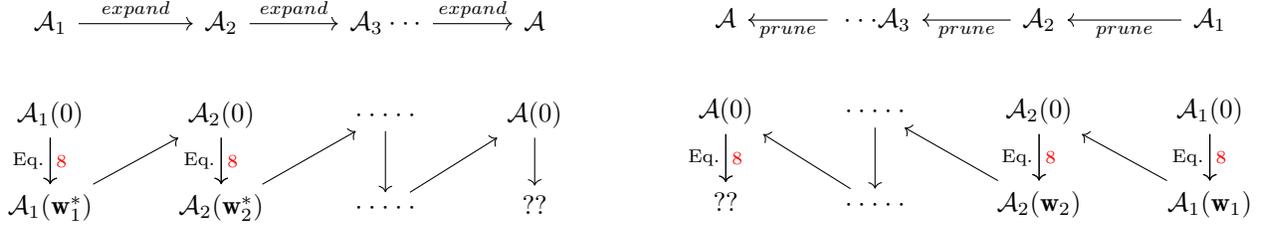
\begin{figure}
\centering
\begin{minipage}{0.48\textwidth}
\centering
\begin{tikzcd}
\CN_1 \arrow[r,"expand"] & \CN_2  \arrow[r,"expand"] & \CN_3 \cdot\cdot \,\cdot \arrow[r,"expand"] & \CN\\
\CN_1(0) \arrow[d,"\text{\ref{param_GD}}"]  \arrow[d,"\text{Eq.}", swap] & \CN_2(0) \arrow[d,"\text{Eq.}",swap] \arrow[d,"\text{\ref{param_GD}}"] & \cdot\cdot\cdot\cdot\cdot \arrow[d,""] & \CN(0) \arrow[d,""]\\
\CN_1(\wst_1) \arrow[ru,""]  & \CN_2(\wst_2) \arrow[ru,"", swap]  & \cdot\cdot\cdot\cdot\cdot \arrow[ru,"", swap]  & ?? 
\end{tikzcd}
\end{minipage}%
\hfill
\begin{minipage}{0.48\textwidth}
\centering
\begin{tikzcd}
\CN  & \, \cdot\cdot \cdot  \CN_3  \arrow[l,"prune"]  & \CN_2 \arrow[l,"prune"] & \arrow[l,"prune"]  \CN_1 \\
\CN(0) \arrow[d,"\text{\ref{param_GD}}"]  \arrow[d,"\text{Eq.}", swap] & \cdot\cdot\cdot\cdot\cdot \arrow[d,""] & \CN_2(0) \arrow[d,"\text{Eq.}",swap] \arrow[d,"\text{\ref{param_GD}}"]  & \CN_1(0) \arrow[d,"\text{Eq.}",swap] \arrow[d,"\text{\ref{param_GD}}"] \\
??  & \cdot\cdot\cdot\cdot\cdot \arrow[lu, ""]  & \CN_2(\weightvector_2)  \arrow[lu, ""]   & \CN_1(\weightvector_1) \arrow[lu, ""] 
\end{tikzcd}
\end{minipage}%
\caption{(\textbf{a}) Where do inclusively expanding $\CN_i$ limit to? $\quad$(\textbf{b}) Can we prune $\CN_i$ without sacrificing performance?}
\label{IAE_schema}
\end{figure}

\subsubsection{Expanding architectures to build better models}\label{expand_Arch}

Let us consider an expanding sequence of architectures $\CN_i$ and their equilibria $\wsit$ under their respective parametric losses $\cL_i$. By Lemma \ref{equilibria}, $\nabla \CL[\CN_i(\wsit)] = 0 \implies \CN_i(\wsit) = \Phi$: at each level, we only need to handle $\Theta_i$ being singular and/or $\Theta_i, \nabla \CL[\CN_i(\wi)]$ being mutually orthogonal. As noted before, while SGD tends to the global minimum in infinite time \cite{sgd_global_convg_geman_86, global_min_SGD_Hajek_88}, in finite time we will at best only be within $\CBep$: we still need a schema that can force us to $\Phi$. Let us consider what the density results like the Universal Approximation Theorem \cite{Horn91} give us. 

\begin{remark}
[\cite{Horn91}]\label{density_sequence}
Let $\CN_i:\R^{M_i} \to G$ be a sequence of architectures s.t. $i < j \Longleftrightarrow M_i < M_j$ and $\CN_i = \CN_j|_{\R^{M_i}}$. Further, assume $ G_\infty = \bigcup_{M_i \in \mathbb{N}} G_{M_i}$ is dense in $G$ and $\FF$ is well-behaved w.r.t. $\CL$ and a solution $\Phi$. Let $\CBep^{M_i}$ be non-empty and $\weightvector^{(*)}_i$ be a global minimizer of $\cL_i[\wi]$ for each $i$. Then $\CN_i(\weightvector^{(*)}_i) \xrightarrow[i \to \infty]{} \Phi$.
\end{remark}

Unsurprisingly, the chain of best models from inclusively larger architectures leads to $\Phi$. However, since empirical optimization flows run for finite time, we need to show that a path through (or near) the first locally critical $\CN_i(\wsit) \in \CBep$ w.r.t. $\cL_i$ can converge to $\Phi$ too. We need to do this constructively by showing a method of extending dynamics from the first well-behaved equilibrium at each $i$ can be enough. Consider a limited example:
\begin{corollary}[Ch. 4, \cite{manyfold_learning}]\label{Linear_IAE}
Assume $\FF$ is strictly well-behaved w.r.t. $\CL$, $\CN_1$ is well-initialised w.r.t. Eq. \ref{param_GD}, and a sequence of linear $\CN_i$ satisfies the conditions of Remark \ref{density_sequence}. Let $\wsiit$ represent the equilibrium under Eq. \ref{param_GD} for initial condition $\weightvector_{{i+1}}(0) = \weightvector^*_{i}$. Then, $\CN_i(\weightvector^*_{i}) \xrightarrow[]{i \to \infty} \Phi$.
\end{corollary}

We need a generalized Corollary \ref{Linear_IAE} that covers nonlinear $\CN$ and merely well-behaved $\FF$.\\

\begin{figure}
\centering

\begin{subfigure}{0.45\textwidth}
\centering
\begin{tikzpicture}[
    neuron/.style={circle, draw, minimum size=18pt},
    ghost/.style={circle, draw, minimum size=18pt, dotted},
    weight/.style={midway, fill=white, inner sep=1pt},
    >=stealth
]

\node[neuron] (x) at (-0.5,0) {};

\node[neuron] (h1) at (3,2.15) {$\sin(\cdot)$};
\node[neuron] (h2) at (3,-2.15) {$\sin(\cdot)$};

\node[ghost, line width=0.75pt] (g1) at (3,4.38) {};

\node at (3,3.5) {$\vdots$};

\node (h3) at (2.95,0)
{$\cNw(\cdot)=w_0+\sum\limits_{i=1}^{a} w_{2i}\sin(w_{2i-1}(\cdot))$};

\node[neuron] (y) at (6.5,0) {};

\node (bias) at (6.5,4.25) {};

\draw[->] (x) -- (h1) node[weight] {$w_3$};
\draw[->] (x) -- (h2) node[weight] {$w_1$};

\draw[->, dotted, line width=0.75pt] (x) -- (g1) node[weight] {{\Huge\textbf{?}}};

\draw[->] (h1) -- (y) node[weight] {$w_4$};
\draw[->] (h2) -- (y) node[weight] {$w_2$};

\draw[->, dotted, line width=0.75pt] (g1) -- (y) node[weight] {\Large$\textbf{0}$};

\draw[->] (bias) -- (y) node[weight] {$w_0$};

\end{tikzpicture}
\end{subfigure}
\hfill
\begin{subfigure}{0.52\textwidth}
\centering
\begin{overpic}[width=\textwidth]{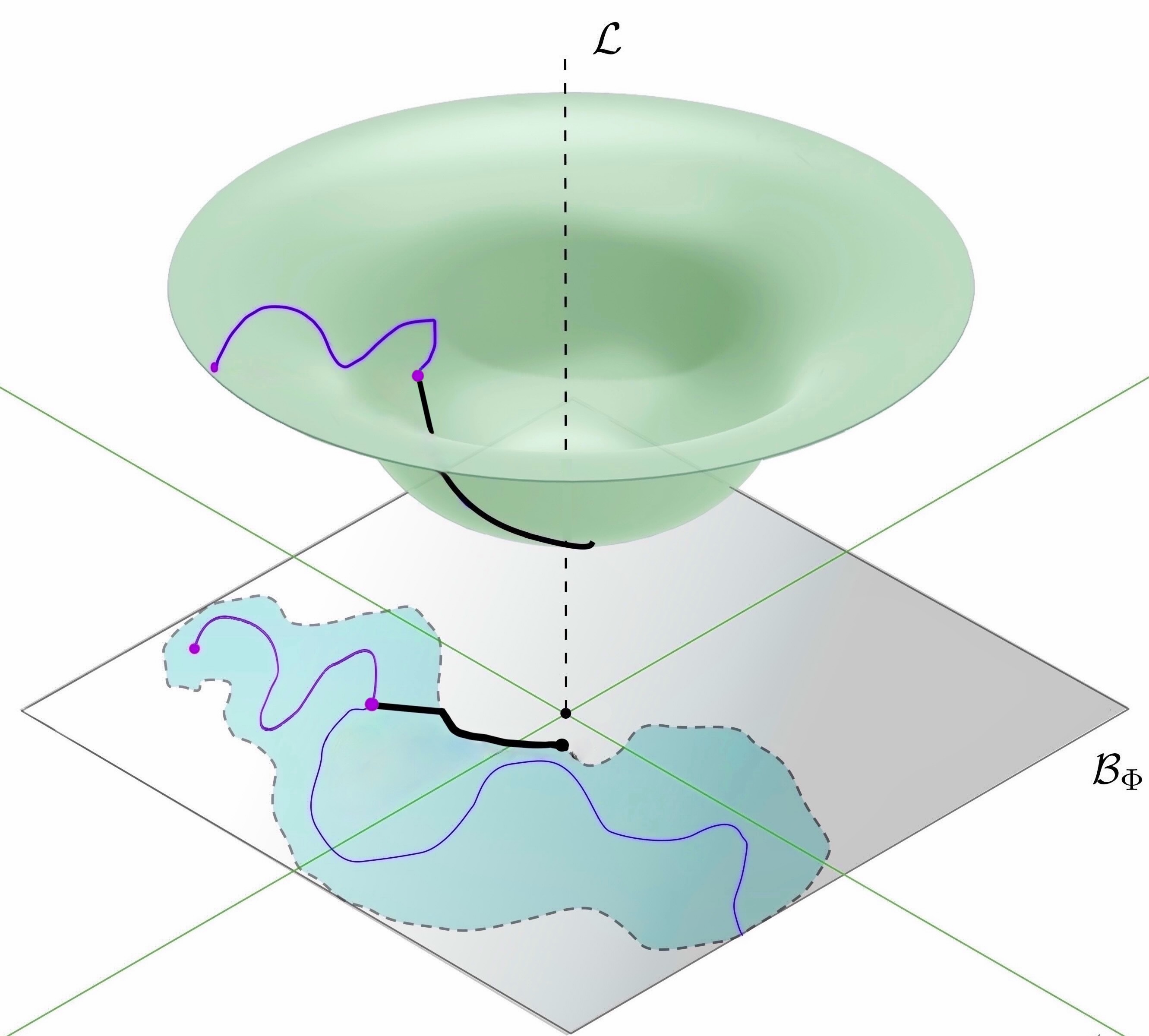}
    \put(32, 25.2){\makebox[0pt][s]{$_{\CN_2(t)}$}}
    \put(56, 30.5){\makebox[0pt][s]{$_{\dot{\CN}_2(t)= -\Theta_2(t)\nabla\CL[\CN_2(t)]}$}}
    \put(64, 34.5){\makebox[0pt][s]{$_{\CN_2(0)=\CN_1(\wst)}$}}
    \put(43.4, 8.5){\makebox[0pt][s]{$_{\leftarrow\,\CBep^2\,\to}$}}
\end{overpic}
\end{subfigure}
\caption{(\textbf{a}) Finite parameter architectures as sub-architectures of larger architectures. Here, $\CN_5(\weightvector') = \CN(\weightvector)$ for all $\weightvector' \in \R^5, \weightvector \in \Rinf$, s.t. $\proj_{\R^5}\weightvector = \weightvector'$ and $w_{2a} = 0$ for $a>2$. (\textbf{b}) Expanded dynamics}
\label{fig:Final_N}
\end{figure}

Consider Ex. \ref{stock2} again: $\cNw(\cdot) = w_0 + \sum_{i = 1}^{a} w_{2i}\sin \Big(w_{2i-1}(\cdot) \Big)$, with $M = 2a+1$. Any such $M-$parameter $\CN$ can be viewed as a sub-architecture of a larger architecture, with the excess parameters set to 0 (alternatively, observe that $G_M \subset G_N$ in Ex. \ref{stock2} if $M < N$). Indeed, $G_\infty$ in Ex. \ref{stock2} is dense in $L^2([a, b])$. This alone doesn't guarantee we can reach models arbitrarily close to $\Phi$ as computers would still not be able to access $\R^\infty$ or $l^2$. But chaining larger architectures in a nested manner could clearly be a solution (Fig. \ref{fig:Final_N}\textbf{b}). Let us formalise this idea:

\begin{definition}[Inclusive Architecture Expansion (IAE)]\label{defn:IAE}
A sequence of architectures $\CN_i \in \cC^2_{pw}(\R^{M_i}, G)$ is called an IAE, if for all $ i, j, M_i, M_j \in \mathbb{N}$: \textbf{(i)} $i < j \implies M_i < M_j$ and $\CN_j|_{\R^{M_i}} = \CN_i$,

\noindent \textbf{(ii)} $i<j \implies \forall \weightvector_i \in \R^{M_i}, \exists \weightvector_j \in \R^{M_j}\text{, } \text{ s.t. } \proj_{\R^{M_i}}(\weightvector_j) = \weightvector_i \text{, }\textbf{ } \mu_j(\weightvector_j) \neq 0$, and $\CN_i(\weightvector_i) = \CN_j(\weightvector_j) = \CN_j(\weightvector_i)$,

\noindent \textbf{(iii)} The direct limit\footref{N_extends} $\GN = \varinjlim  G_{M_i} = \bigcup_{M_i \in \mathbb{N}}G_{M_i}$ is a subspace of $G$. 
\end{definition}

IAEs formalize how we would intuitively expand architectures while respecting their “structure'' and preserving the performance of our existing model. Condition \textbf{(ii)} implies that for a sequence of architectures to qualify as an IAE, it has to be able to open new gradient directions for optimization without changing the model, even as the architecture is changed as new parameters are added (in a potentially randomly manner). 

Note that $\CN_j(\weightvector_i)$ represents $\CN_j$ applied to the canonical extension of $\weightvector_i$ to $\R^{M_j}$.
These assumptions are satisfied by a large class of practically useful architectures. If $\CN_i$ are nonlinear, \textit{\textbf{(ii)}} can also be fulfilled by choices $\weightvector_j \neq \weightvector_i \in \mbbR^{M_j}$, but $\proj_{\R^{M_i}}(\weightvector_j) = \weightvector_i$. Could we thus speak of a final architecture that accepts infinite parameters? 

\begin{lemma}[Final architecture, Ch. 4, \cite{manyfold_learning}]\label{final_arch}
Let $\CN_i$ be an IAE. Then there exists a unique $ \CN \in \cC^2_{pw}(\Rinf, G)$ s.t. for all $i$, for all $\weightvector \in \R^{M_i}$, $\CN(\weightvector) = \CN_i(\weightvector) $. 
\end{lemma}
\begin{remark}[Abuse of notation]
    Since $\CN_i$ for all $i$ in an IAE may be viewed as sub-architectures within some final architecture $\CN$ and $i < j \implies G_{M_i} \subset G_{M_j}$, we may drop the first $i$ from $\CN_i(\weightvector_i)$, whenever $\CN(\weightvector_i)$ makes the architecture level clear. Similarly, $\cL_i(\wi) \to \cL(\wi)$, $\Theta_i(\weightvector_i) \to \Theta(\weightvector_i)$, etc. Whenever relevant, a distinction will of course be made between $\wi \in \R^{M_i}$ and $\wi \in \Rinf$. Finally, if we are analysing only a single fixed architecture $\CN$, we can always view it as the final architecture of the directed system $\CN_i = \CN, \forall i \in \N$.
\end{remark}

As an example, assume an $M$ parameter Ex. \ref{stock2} solving Ex. \ref{stock} limits to a local critical point of Eq. \ref{param_GD} (implying $\mu(\weightvector) \to 0$ or $\nabla \CLN \in \ker({\Theta(\wst)})$, since $\nabla \CL[\CNwst]$ is 0 only at $\CNwst = \Phi$). Let us expand this by adding $2$ new parameters s.t. $w_{2a+2} = 0$ but $w_{2a+1} \neq 0$, thus ensuring that $\mu(\weightvector_{M+2}(0)) \neq 0$. The existence of such new parameters, that can open new optimization directions as a guarantee, while extending from old parameters, means we have a computationally viable way to take the limit $M \to \infty$. Indeed, this trick is generalisable and we may add as many extra parameters as wanted to any $\CN$ without losing the original model: if all new linearly dependent parameters (say $l_i)$ are set to 0 at expansion, while all other parameters are chosen to ensure non-zero $l_i$ derivatives.

\subsubsection{A prototype for a universal convergence theorem}
Further, $G_{M_i}$ usually approaches a dense linear subspace structure as $M_i \to \infty$ (\cite[Lemma 4.9]{manyfold_learning}): the \textit{manifolds} $G_{M_i}$ flatten out as $G_\infty$ is approached. In this regard, our schema is improving upon the classical NTK approach, by building up to the infinite parameter limit as optimisation is carried out, rather than optimising in the limit.

This point underpins the linearized/convexified optimization regime seen in the kernel training setting (convexification being the result of stronger conditions on $\CL$ than LI). A similar phenomenon may allow us to approach the global minimum, as sub-optimal critical points are destroyed by the IAE \cite{Loja_to_ML_2022}. Let us define the schema by which dynamics may be restarted with the guarantees that new viable gradient flow directions will be found:

\begin{definition}[Error Correcting Inclusions (ECI)]\label{ECI_defn}
    An IAE $\CN_i$ associated with the following sequence of dynamical systems is called an ECI:
    \begin{equation}\label{CGD_central}
    \begin{aligned}
    & \frac{d}{dt}\weightvector_{i+1}(t) = -\nabla\cL[\weightvector_{i+1} (t)] \implies \frac{d }{dt}{\CNt} = -\Theta_{i+1}(t) \nabla \CL[\CNt] \\
    & \text{with } \weightvector_{i+1}(0)\quad\text{s.t. }\proj_{R^{M_i}}(\weightvector_{i+1}(0)) = \wst_{i}, \textbf{ } \CN(\weightvector_{i+1}(0)) = \CN(\wst_{i}) \\
    \end{aligned}   
    \end{equation}
\end{definition}

\begin{theorem}[Ch. 4, \cite{manyfold_learning}]\label{vtk_central_2}
    Assume $\FF$ is well-behaved w.r.t. $\CL$, $\CN$ is the final architecture of an ECI $\CN_i$, $G_\infty$ is a dense subspace of $G$, $\CN(\weightvector_1(0))$ is well-initialized, and $\cL$ is locally well-behaved at each $\weightvector^*_i$. Then $\CN(\wst_i) \xrightarrow{i \to \infty} \Phi$.
\end{theorem}
$\CN(\weightvector_1(0))$ being well-initialised is the critical assumption above: as noted before, most $\FF$ are easily shown to be well-behaved \cite[Ch. 3]{manyfold_learning} and useful architecture choices usually do satisfy the density and expansion conditions (entirety of Table \ref{fig:method_and_summary} does). Further, as most practical architectures are piecewise analytic over bounded domains, each expansion generates locally well-behaved critical points (Corollary \ref{example-Li}). Note that regions of discontinuity represent a different kind of difficulty than criticality, but modern sub-gradient methods handle it well \cite{FoCM_20_stoch_sub_Tame_Loja_Damek}.

We posit that Thm. \ref{vtk_central_2} is a prototype for a universal convergence result: it blends the solvability of a problem with universal approximation capabilities of a method to guarantee arbitrarily good solution models. Qualitatively, it says that if we keep expanding the model sets $\Gi$ at local solutions $\CN(\wsit)$, we can construct a path limiting to $\Phi$, even if $\Phi \neq G_M$ for any $M \in \N$. Crucially, we do not begin learning in the over-parametrised regime \cite{jacot_ntk_18, ntk19jae}.

Fig. \ref{fig:Final_N}\textbf{b} presents a visual summary of how ECIs lead us to target solutions. Note that ECIs are used \textbf{after} SGD or some other coarse-graining technique has produced a good model in $\CBep$: it does not have to deal with the geometry of $\CL$ in the entirety of $G$, unlike if we started with a larger architecture in the first place. Further, the expected performance of a highly parametrised NN may not always change meaningfully with even a doubling in size - a lot still depends on the chosen problem/architecture combination. Like conventional NTK theory, ECIs only guarantee consistently better performance on $\CL,\cL$ and eventual convergence - but in a feasible manner.

\subsubsection{Sparsifying or pruning architectures to retain or improve performance at lower costs}\label{prune}

We have learnt how to expand $\CN$ without losing the model. We now sketch two sets of ideas that allow us to prune $\CN$ without losing performance: \textbf{(i)} a reverse ECI process that prunes away parameters whenever the model manifold becomes degenerate \textbf{without} interrupting the model flow (Fig. \ref{IAE_schema}\textbf{b}) and \textbf{(ii)} statistical viewpoints on the parameters comprising $\weightvector$ that eliminate them \textit{en-masse} post a certain amount of learning \cite{fra19, mae21, reddog22}. The focus is brief because these ideas are more heuristic than rigorous and the framework is quite limited on them as of now.

In $(\textbf{i})$, we wish to move the optimisation dynamics to a smaller model set without losing performance \textit{in real time}. We can then invest the saved costs in ECIs if needed. Note how ECIs expand $G_M$ at $\CNwst$ to restart paused dynamics. How do we identify $\cNw$ that sparsify the geometry of $G_M$ without pausing the optimisation dynamics?

The same nonlinearities that complicate analysis also allows degeneracies we can use to scale $\CN$ up or down. Consider Example \ref{stock2} at some non-critical $\wt$ such that $w_4(t) = 0 \implies \frac{\partial \CN}{\partial w_3 } = 0$. If we eliminate $w_3, w_4$ for future optimisation from time $t$, we do not lose the model and optimise in the lower cost $
\R^3$ regime. When optimisation pauses due to a critical point $\wst \in \R^3$, we can use an ECI to restart the flow in $\R^5$. This motivates:
\begin{definition}[\textit{in situ} Pruning \cite{dog23_3}]
    Let $\weightvector(t), \CN(t)$ represent parametric and model optimisation flows and $\weightvector_j \subsetneq 
    \weightvector(t)$ represent the flow for a proper sub-collection of parameters. We call $\CN_{j}$ and $\weightvector_j(t)$ an \textit{in situ} pruning of $\CN$ at time $t$ if: (\textbf{i}) $\CN_{j} = \CN|_{\R^{M_j}}$ and (\textbf{ii}) $\weightvector(t)$ and $\weightvector_j(t)$ are not critical points for $\cL, \cL_j$ respectively.
\end{definition}
\begin{remark}
    Let $\CN$ be an architecture that is linearly dependent on some parameter $w_l$. Let $\weightvector^{(l)}$ represent the parameters that are exclusively within the argument of $w_l$. Let $\R^{M_j}$ represent the space of parameters not containing $w_l$ and $\weightvector^{(l)}$. We can prune away $w_l$ and $\weightvector^{(l)}$ \textit{in situ} every time $w_l(t) = 0$ and $\nabla \CL[\CN(\weightvector_j)]$ is not orthogonal to $\Theta_j$.
\end{remark}

Statistical viewpoints focus on identifying sparse sub-networks in deep NNs that provide essentially similar performance when trained in isolation \cite{fra19}. Indeed, those identified over one problem-method combination can often be used in others \cite{meh19}. For predominant methods like Iterative Magnitude Pruning (IMP) \cite{tanaka2020pruning}, Renormalisation Group (RG \cite{wil71A}) theory provides a useful foundation for considering these aspects of pruning \cite{manyfold_learning, abu_23, red22iclr, reddog22}. Since we did not develop novel ideas or tools in RG theory, we will only comment on how we can leverage it for pruning.

Let $\CN$ be an architecture over a set of activations $\bm \phi$ and parameters $\weightvector$, $\cL[\wst, \bm{\phi}]$ be the post training performance of $\CN$ on a problem-method combination, and $\delta_i = M_i/M$ be the proportion of left-over parameters after $i$ pruning cycles on $\CN$. IMP variants work by: (\textbf{i}) identifying and eliminating \textit{dead} parameters/activations (essentially 0-valued and/or static) at the end of training, (\textbf{ii}) rewinding survivors to their pre-training values and training again. 

We can setup the IMP algorithm as an operator and view $\cL_i[\wsit, \bm\phi_i] = \mathcal{I}^i \cL[\wst, \bm{\phi}]$. Empirical observations show that $ \cL_c[\wst_c, \bm{\phi}_c]  \approx \cL[\wst, \bm{\phi}] $ up to a critical $\delta_c$ and follows a power law after: $\cL[\wsit, \bm \phi_i] 
\sim (\delta_c - \delta_i)^{-\gamma}$ \cite{fra19}. 

In \cite{abu_23, reddog22}, $\mathcal{I}$ was rigorously shown to be an RG process. The spectrum of this operator identifies \textit{relevant} and \textit{irrelevant} parameters by invoking existing tools from RG theory - a fact of limited value, since $\CN$ has to be fully trained at least once to allow this identification. However, RG theory also allows us to build universality classes of problem-method combinations, which allows for profitable transfer of the identified sub-networks across tasks.

\section{Limitations and Future Directions}\label{limitations}

We have built a framework that covers a large universe of problems and methods \cite[Ch. 3]{manyfold_learning}. However, it would benefit from deeper results and more engagement in the following directions (in increasing order of difficulty):
\begin{itemize}
    \item Target spaces: weakening the requirements on $G$ to it being a metric and perhaps merely topological space. 
    
    \item Model sets: allow more geometrically complex model-sets within $G$ than essentially differentiable manifolds.

    \item Costs: theory can only be directly tested up to $M \sim 10^4$ and the tests cost more than the actual optimisation.

    \item Update rules: incorporate stochastic $\wt, \CNt$ dynamics alongside essentially differentiable ones. Further, while methods with $\CL$ as a Lyapunov map are gradient flows, it would help a lot more if we also knew how.

    \item Existing theories: other ideas in statistical and deep learning theory can be built atop the framework but it is not immediately clear what we could add to those ideas beyond the facets we have discussed already.

    \item Well-initialization: powerful assumption and we offer nothing beyond the “expectation'' SGD can provide it.

    \item Data: acknowledged to have a strong implicit role (setting up $\CL$, choosing $\CN$, etc) that we offer no insight on. 
\end{itemize}
We ask $G$ and $\CN$ to have differentiable structure partly because it is useful in proving the non-degeneracy of the model sets. However, gradient flow theories exist and are well understood for metric spaces without Riemannian structure \cite{Savare_GD_05}. If $G$ could be a metric (or perhaps even topological) space without losing our results on $G_M$, we would directly access a wider span of problems. Similarly, if $\CN$ could parametrise $G$ using more complex structures, without losing the geometrical and dynamical properties we have relied upon for convergence results, we would significantly expand the validity and usefulness of these ideas.

Indeed, even though $\mu(\weightvector)$ is estimable at each $\weightvector$ (Sec. \ref{arch_features}), the cost of such estimation can scale as $\mathcal{O}(M^3)$ at each $t$. A direct approach to such estimation can rapidly become unviable (modern ML architectures can have $M > 10^{12}$). It is perhaps unsurprising that we have only trialled and built applications in the regime of low-cost scientific modelling, where $M \sim 10^2$ (SHAPER \cite{dog23}, FINDER \cite{trajan25}) to $10^4$ (NN DE solvers \cite{dog23_3, weinan18, mar20, Sirig18, shin_karni_20}, small transformers \cite{reddog22}, and tracking optimization dynamics \cite{diet20, dog20, dogred20}). In this setting, our results can be actually implemented to guide choices and heuristics. If $M \gtrsim 10^5$, hope and indirectly probe are the only concrete choices. Further, we only have local convergence rates for a fixed $i$ in an ECI: convergence rates to $\Phi$ w.r.t. $i$ are unknown.

The framework is built with piece-wise regular (essentially differentiable) flows over $\R^M, G_M$ in mind. However, while it can directly support stochastic and even more irregular optimisation dynamics, we have no insight on how or why it would be beneficial to do so from within this perspective. Similarly, while there is a coordinate chart under which \textit{good} learning dynamics are  gradient flows, unless we explicitly obtain that chart, the practical value of Remark \ref{loss_para} is limited to loose bounds. Indeed, if $\dim(G) < \infty$, most of our results are inferior to existing ones \cite{Convex_Opt_04}.

Stochastic gradient descent (SGD) offers us a weak guarantee: we can \textit{expect} to achieve well-initialised states in \textit{due time} (convergence in law). However, we have no generic ways of validating if that actually happens: the desired result would be to say that well-initialisation happens almost surely for any \textit{reasonable} set of initial parameters. Applications can often involve setting up rough estimates using SGD or other external means and then optimizing from them. For example, in \cite{dog23}, optimization over real datasets happens after obtaining rough pre-estimates using cheaper methods. Similarly, in \cite{dog23_3}, ECIs are activated at low loss values \textit{hoping} that implies $\CBep$ has been reached. Clearly, this is not a completely crippling issue - but it is a critical gap and we have offered nothing on it. 

To make matters worse, we have implicitly assumed that more practical matters like conditioning, regularisation, generalisation, etc, are being taken care of by choosing or transforming $\FF, G, \CL, \CN, \CX$, etc, ``appropriately''. This is a significant assumption, and often \textbf{the} non-trivial issue in a practical problem. The available data is usually our strongest tool in making these choices. A better unified theory would have a component that caters to how such decisions can be made in a data-driven manner. For example, $\FF$ may be too ill-conditioned in its original form or $G$ too sensitive a space to initial conditions. A sequence of data-driven transforms on $\CL$ or $\CN$ may exist that converts  them from being computationally intractable to feasible: this perspective gives us few hints on what those may be. As of now, other ideas in statistical and deep learning theory are far more useful in such considerations \cite{rob21, stat_learn_theory_vapnik00}.

Our distinctive strategies lay in separating the problem and the method at the outset and then showing how their combinations result in properties that can be quantified explicitly and \textit{a priori}, in terms of their independent ones. Future work will focus on resolving the challenges we have noted as rising out of such choices and developing practical applications based on the strengths of our perspective.

\section{Conclusions and related works}

Our goal was to sketch how a unified theory of learning could be built. We assembled a prototype by (\textbf{i}) establishing when a problem is solvable in principle, (\textbf{ii}) encapsulating the balance between linear and nonlinear methods to solve it, and (\textbf{iii}) investigating how the problem and method mix to govern parametric optimisation profiles. 

Unsurprisingly, our endeavours are situated well within the span of a rich existing literature, especially of the last two decades. This work is a synthesis of existing ideas on convex optimization \cite{Convex_Opt_04}, NTK theory \cite{ntk_imp_19, jacot_ntk_18, ntk_imp_18}, its extensions beyond least square losses \cite{cohen23} and infinite widths \cite{seleznova22b, seleznova22a, Inf_Wid_useless_23}, \textit{a priori and a posteriori} error bounds for NN DE solvers \cite{dog23_3, Posteriori_23, Deeponet_error, FNOs_20, mar20}, characterization of optimization dynamics \cite{townsend_23, hoop_stuart23, diet20, dogred20,  ntk_vs_dl_20, ntk19jae, exp_convg_schwab_23}, constraint enforcement in scientific ML \cite{dog23_3, mar20}, deep learning failure regimes \cite{Mahoney_Pinn_fail_21, Mariia_NTK_Collapse_23, sifanwang20}, etc. Consequently, our claim and aim is not some radical breaking of the mould: it is the construction of a broad framework that seamlessly encompasses existing empirical observations and theoretical findings, while unearthing new methods and questions.

We broke down the utility of both nonlinear and linear learning methods into three types of properties: density of $G_M$ in the target space $G$, geometrical and dynamical properties inherited by $\cL$ and $G_M$ from $\CL$ and $\CN$, and the optimisation dynamics used to force the models to $\Phi$. 
We did not establish novel density results for model sets since Universal Approximation Theorems exist for many choices of $G$ \cite{Cyb89, GKP20, UAP_Mishra_21, Deeponet_error,  UAP_Bochner_space_Schmocker_23, Bouned_Relu_SHEN_22}. However, they were augmented by invoking manifold like structures on $G_M$ to better study the possible dynamics.

Solvable (or well-behaved) problems were shown to transmit convergence properties to their parametrised analogues by adaptations on the Lojasiewicz Inequality (LI, \cite{Loja_og_63}), especially when paired with the manifold structures invoked on  the model sets. The results on forcing a path to the solution by going through a sequence of expansions on $\CN$ adapt some conventional ideas from Neural Tangent Kernel (NTK \cite{jacot_ntk_18, ntk19jae}) and over-parametrised regimes. However, our focus was on the parameters and computational feasibility, while bypassing an explicit role for data.

Thm. \ref{vtk_central} and \ref{vtk_central_2} are key distinguishing features, since the limited setting of over-parameterized architectures has been previously used to characterize optimization dynamics and guarantee arbitrarily good solutions.  However, they  struggle to describe empirical optimization dynamics, being fundamentally rooted in a computationally inaccessible setting \cite{Mariia_NTK_Collapse_23} that forces an essentially convexified framework at odds with ground realities. 

Recent modifications \cite{Loja_to_ML_2022, seleznova22b, seleznova22a} have gone further to address these limitations: indeed, we are strongly influenced by \cite{Loja_to_ML_2022, LI_metric_proximal_grad_Lorenzo_24}. However, they are either limited to problems over finite datasets \cite{seleznova22b, seleznova22a}, still rooted in over-parameterized regimes for some results \cite{Loja_to_ML_2022}, or concerned with establishing results only in the parametric setting independent of the available results or structure on the original problem \cite{LI_metric_proximal_grad_Lorenzo_24}. Similarly, while the idea of learning as a dynamical system is not new \cite{diet20, dogred20}, especially in convex settings \cite{Convex_Opt_04, helmke_moore_1994_optimization}, we are considering optimization in the versatile setting of nonlinearly generated model sets in possibly infinite dimensional manifolds. We are agnostic to the number of parameters, show how to scale architectures up and down without losing performance, and approach the over-parametrised regime through a feasible and self-contained schema.

\section*{Acknowledgements}

This work is supported by the National Science Foundation under Cooperative Agreement PHY-2019786 (The NSF AI Institute for Artificial Intelligence and Fundamental Interactions, http://iaifi.org/). This research was also funded by the President’s
PhD Scholarships at Imperial College London and by the EPSRC Centre for Doctoral Training in Mathematics of Random Systems:
Analysis, Modelling and Simulation (EP/S023925/1). 

This work grew cumulatively from the lecture notes prepared for the following events and/or host institutions (organisers/participants at these events have our strongest gratitude for their kind interest and passionate discussions):

1. Computational and Applied Mathematics Laboratory at ETH Zurich (May 2023)

2. Applied and Computational Analysis Seminar at the DAMTP, University of Cambridge (March 2024)

3. NSF IAIFI Colloquium, Department of Physics, MIT (April 2025)

4. Advanced Summer School on Quantum Field Theory and Quantum Gravity, CERN/ICISE (June 2025)

We are also thankful to Mr. Evangelos Kotzafilios (Imperial College London), Prof. Jeroen S. W. Lamb (Imperial College London), Dr. Martin Peev (University of Oxford), Dr. William T. Redman (Johns Hopkins), Dr. Kevin N. Webster (Imperial College London), and Dr. Dominic Wynter (UT Austin) for their insightful discussions.

\renewcommand\bibname{References}
\bibliography{References}
\bibliographystyle{plain}

\appendix \label{math_appendix}

\newpage

\section{NOMENCLATURE}\label{nomenclature}
\subsection{Fundamental objects of interest}

\renewcommand{\arraystretch}{1.45}
\begin{tabular}{ll}
\textbf{Symbol} & \textbf{Description} \\ \hline

$\FF$ & Problem of interest: maps between complete separable metric spaces $G, H$. \\

$G$ & Target Space: space to which the possible solutions to the problem of interest belong.\\

$\CL$ & Nominal Loss: maps from $G$ to $\R$ with their minimum at the true solution(s) $\Phi$. \\

$\Phi$ & Solution(s): objects in $G$ to be modeled. Often themselves maps between other spaces. \\

$\weightvector$ & Parameters: trainable elements in spaces $\R^M, M \in \mathbb{N} \cup \infty$, to be optimized. \\

$\CN$ & Architectures: maps from {parameters} to model elements in $G$. \\

$\cNw$ & Models: elements in $G$ that are also images of a parametrized architecture. \\

$G_M$ & Model Set: Set of models generated by $\CN$, given as
$G_M \eqdef {\CN(\weightvector): \text{ $\weightvector \in \mbbR^M$}}$. \\

$\cL$ & Parametric Loss: maps parameterizing $\CL$ over $\R^M$ through $\CN$. \\

$\weightvector^*$ or $\wst_M$ & Limit point of the dynamical flow $\wt \in \R^M$, if it exists. \\

$\weightvector_M^{(*)}$ & Global minimum of $\cL$ over $\R^M$. \\

$\Theta$ & Generalized Neural Tangent Kernel for $\CN$. \\

$\theta$ & Shares its non-zero spectrum with $\Theta$. \\

$l^2$ & Space of square-summable real sequences. \\

$\Rinf$ & Inductive limit of convergence (Direct limit) of the direct system $\{\R^M\}_{M \in \N}$. \\

\end{tabular}

\subsection{Notation}\label{notate}

\renewcommand{\arraystretch}{1.45}
\begin{tabular}{ll}
\textbf{Symbol} & \textbf{Description} \\ \hline

$D A$ & Fréchet derivative of a map $A$ over an arbitrary space $G$. \\

$\mrmD A$ & Fréchet derivative of a map $A$ over real space $\R^M$. \\

$\nabla A$ & Gradient of a map $A$. \\

$\nabla^2 A$ & Hessian $\nabla(\nabla A)$ of the map $A$. \\

$\Wkp(X)$ & Sobolev space of order $k$ under $p$ norm, where $X \subset \R^N$ is an open Lipschitz domain. \\

$A^\dagger$ & Adjoint of a linear operator $A$. \\

$\sigma \Big(A\Big)$ & Spectrum of a linear operator $A$. \\

$\sigma_+ \Big(A\Big)$ & Singular values of a linear operator $A$. \\

$\Ndw$ & Derivative $\mrmD \CN$ of an architecture $\CN$ w.r.t. its parameters. Its adjoint is $\Ndwdag$. \\

$\CNww$ & Hessian of $\CN$ w.r.t. its parameters. \\

$\cC^k(X)$ & Maps with locally continuous derivatives up to order $k$. \\

$\cC^k_{pw}(X)$ & Maps with piecewise locally continuous derivatives, with a finite number of discontinuities  over $X$. \\

$\overline{X}$ & Closure of $X$. \\

$\varinjlim A_i$ & Direct limit of a family of indexed objects $A_i$ (usually algebraic objects). \\

$f|_X $ & A map $f$ restricted to a sub-domain $X$ \\

\end{tabular}

\subsection{Abbreviations}\label{Abbreviate}

\renewcommand{\arraystretch}{1.3}
\begin{tabularx}{\textwidth}{l X l X}
\textbf{Symbol} & \textbf{Description} & \textbf{Symbol} & \textbf{Description} \\ \hline
i.e. & that is & Ex. & Example \\
w.r.t. & with respect to & Defn. & Definition \\
s.t. & such that & Ch. & Chapter \\
a.e. & almost everywhere & Fig. & Figure \\
p.s.d. & positive semi definite & Sec. & Section \\
w.l.o.g. & without loss of generality & Eq. & Equation \\
IAE & Inclusive Architecture Expansions & Thm. & Theorem \\
ECI & Error Correcting IAE & App. & Appendix \\
\end{tabularx}

\section{Supplementary Remarks}\label{pert_note}
\subsection{Universal Approximation Property (UAP)}\label{UAP}
The existence of ``good'' models is guaranteed for a large class of solution(s)/target objects by variants of universal approximation results. In particular, let $\activation$ be a non-polynomial, Lipschitz, and a ``sufficiently regular'' nonlinearity (usually we need $\activation \in \cC^2_{pw}(\R)$ or less). A sequence of works has established that fundamental units like $\activation$ in combination with a countable number of parameters generate a set of Neural Network models that can be dense in the space of $\cC^k$ functions \cite{Cyb89}, in Sobolev spaces \cite{GKP20}, in Wasserstein \cite{UAP_Wasserstein_space_Lu_20} and Bochner \cite{UAP_Bochner_space_Schmocker_23} spaces, and in the space of continuous nonlinear operators between Banach spaces \cite{Deeponet_error, UAP_Mishra_21}. Indeed, arbitrarily good models sometimes ``exist'' even when $\CN$ is allowed only a finite number of parameters \cite{bounded_ab}. However, finding them is another matter\footref{well_init}.

\subsection{Well Initialization}\label{well_init}
Definitive answers on how to practically obtain ``well-initialised'' models or validate such a label are beyond the scope of this work. In upcoming work, we will aim to bound the search times expected under stochastic gradients (\ref{SGD_uses_lin_vs_nlin_N}) to reach neighbourhoods of solutions from arbitrary initial conditions, under some structural assumptions on the loss functions being used to guide the dynamics. In practice, ``well-initialization`` is usually achieved by \textit{pre-training} on limited empirical data \cite{dog23}, using parameters of models trained on similar tasks (\textit{transfer learning} \cite{red22iclr}), hoping that the convergence in probability results on stochastic descent allows good models within feasible costs (App. \ref{SGD_uses_lin_vs_nlin_N} and \cite{sgd_global_convg_geman_86}), etc. Validation of such models is a nascent discipline \cite{comp_assist_2019}.

\subsection{UAP is countability}\label{UAP_sep}

When $G, H$ are Hilbert spaces (or their tangent spaces are), requiring a countable Hilbert dimension is equivalent to assuming separability: both assumptions are not generically true even for Hilbert spaces. However, they are implicit if we are using discrete methods. Absent {separability}, there would be elements in $G$ that can't be computationally modelled to arbitrary accuracy.

Since we already assumed $G$ usually has a countable Hilbert dimension, this assumption essentially requires that the union over the space of finite models is rich enough to be dense in $G$. It is a universal approximation property (UAP) assumption, which is reasonable since UAP for NNs has been proven for large classes of function spaces and NNs \cite{han17, Kidger19, yar17, yar18}. Indeed, in many ways, such results are one of the bedrocks of machine learning.

\newpage
\section{{Fundamental Definitions and Lemmas}}\label{fund_lemmas}
The definitions and results in this section are standard in their respective fields (or direct implications thereof).  We begin by making precise the notion of a manifold modelled on Banach spaces:

\begin{definition}[Banach Manifolds]\label{Banach_manifold}
Let $G$ be a Hausdorff, second countable topological space and $k \in \mathbb{N} \cup \{\infty\}$. We say $G$ is a separable $\cC^k$ Banach manifold with an atlas $\{(U_\alpha, \varphi_\alpha)\}_{\alpha \in A}$, where each $U_\alpha \subset G$ is open and $\varphi_\alpha : U_\alpha \to \varphi_\alpha(U_\alpha) \subset X_\alpha$ is a homeomorphism onto an open subset of some Banach space $X_\alpha$, if:
\begin{enumerate}
    \item The charts cover $G$: $\bigcup_{\alpha \in A} U_\alpha = G$.
    \item For overlapping charts $(U_\alpha, \varphi_\alpha), (U_\beta, \varphi_\beta)$, the transition map $\varphi_\beta \circ \varphi_\alpha^{-1} : \varphi_\alpha(U_\alpha \cap U_\beta) \to \varphi_\beta(U_\alpha \cap U_\beta)$ is a $C^k$-diffeomorphism
\end{enumerate}
$G$ is a connected manifold if for some fixed Banach space $X$, we may set $X_\alpha = X$ for all $\alpha \in A$. We say $G$ has a strong Riemannian structure or admits a strong Riemannian metric or is a Banach-Riemannian manifold if $X_\alpha$ are reflexive spaces with a smoothly varying isomorphism $r_x: X_\alpha \to (X_\alpha)^*$. Finally, we say $G$ is a Hilbert manifold if $X_\alpha$ are Hilbert spaces. We will denote the tangent space of $G$ at $g$ by $T_g G$.
\end{definition}
\begin{remark}
    Given a Banach manifold with a strong Riemannian structure, we can make the identification $\nabla \CL[g] = r_g^{-1}(d\CL[g])$ to obtain a precisely defined gradient flow for $\CL$. Hilbert manifolds generate this identification canonically via the inner product. Separable infinite dimensional metric Banach manifolds can be embedded as open subsets in a Hilbert space \cite{Banach_manifold_embed_Hilbert_spaces_HENDERSON_70}.
\end{remark}

\begin{remark}
    Note that if $G$ is a Hilbert manifold, Riesz representation generates a canonical identification between $D\CL[g]$ and $\nabla \CL[g]$, in turn allowing the kinds of canonical gradient flows given in Eq. \ref{idealized_G_flow}. However, if $G$ only admits a Banach–Riemannian or weaker structure, we need to pick the “right'' duality or $r_g$, since there is no canonical isomorphism we can use to transform Eq. \ref{generalised_Banach_LI} into Eq. \ref{generalized_gradient_LI}. $G$ will usually be a Hilbert space, unless stated otherwise.
\end{remark}
\vspace{\baselineskip}

\subsection{Fundamental concepts in dynamical systems theory}
\begin{definition}[$\omega$ limit set]\label{omega_limit_sets}
Let $G$ be a metric space and $g:[0, \infty) \to G$ represent a bounded ODE solution to:
\[
\frac{d}{dt}g + F[g] = 0, \qquad\qquad F \in \cC(X, G), \qquad\qquad X \subset G\text{ is open and connected}
\]
We define its omega-limit set $\omega(g) \subseteq \{u: F[u] = 0\}$ as: 
\[
\omega(g) = \{g^* \in \R^M: \exists t_n \to \infty \text{ s.t. } g(t_n) \to g^* \} 
\]
\end{definition}

\begin{definition}[Gradient like system]\label{grad_like_flow}
Let $G$ be a Hilbert space, $\CF \in \cC(G, G)$, $\text{Cr}(\CF) = \{u \in G: \CF[u] = 0\}$, and consider the following ODE:
\begin{equation}\label{gradient-like}
\frac{d}{dt}g(t) + \CF[g(t)] = 0    
\end{equation}
We call Eq. \ref{gradient-like} a gradient like system w.r.t. a strict Lyapunov function $\CL$, if $\braket{\CF[g], \CL'[g]} > 0$ for all $g \notin \text{Cr}(\CF)$.
\end{definition}

\begin{definition}[Lyapunov function]\label{lyapunov}
Let $\dot{g} + F[g] = 0, F \in C(G, G)$ describe some dynamical system over $G$. We call $\CL \in \cC^1(G, \R)$ a Lyapunov function for this system, if $\braket{\CL'[g], F[g]} \geq 0$ for all $g \notin \{u: F[u] = 0\}$. $\CL$ is a strict Lyapunov function if $\braket{\CL'[g], F[g]} > 0$ for all $g \notin \{u: F[u] = 0\}$.
\end{definition}

\begin{definition}[Immersion]\label{immersions}
A map $F: A \to B$ is an immersion in $B$, if $D_aF:T_aA \to T_{F(a)}B$ exists and is injective for all $a \in A$. The image $F(A)$ is called its immersed submanifold in $B$. We say the set $F(\mathcal{B}_a)$ is locally an immersed submanifold around $F(a)$, if it satisfies the definition on some neighbourhood $\mathcal{B}_a$ of $a$. Note that for a Banach space $A$, $T_a A \simeq A$.
\end{definition}

\subsection{Some fundamental algebraic and topological concepts}

\begin{definition}[LB Spaces]\label{LBspace}
Let $G_i$ be a sequence of Banach spaces s.t. $G_i \subset G_{i+ 1}$ and the inclusion map $i_n: G_i \to G_{i+1}$ is continuous and linear. $G = \cup_{i \in \N} G_i$ is called an LB space if it is imbued with the finest locally convex topology s.t. each inclusion $i_m: G_m \to G$ is continuous.
 
The finest locally convex topology is the requirement that for all locally convex space $H$, a linear map $T: G \to H$ is continuous iff the restriction $T: G_n \to H$ is continuous for each $n$.
\end{definition}

\begin{definition}[Direct Limit]\label{N_extends}
    Let $(I, \leq)$ be a directed set, $\{A_i: i \in I\}$ be a family of objects indexed by $I$, and $f_{ij}: A_i \to A_j$ be a homomorphism for all $i \leq j$ such that $f_{ii}$ is the identity on $A_i$ and $f_{ik} = f_{jk} \circ f_{ij}$ for all $i \leq j \leq k$. Then, $\{ A_i, f_{ij} \}_{i \leq j}$ is called a directed system and the direct limit is the quotient
\[
\varinjlim A_i = \left( \bigsqcup_{i \in I} A_i \right) \Big/ \sim
\]
where $ A_i \ni a_i \sim a_j \in A_j$, if and only if there exists $k \geq i,j \text{ such that } f_{ik}(a_i) = f_{jk}(a_j)$.
\end{definition}

\begin{definition}[Inductive Limit of Convergence on $\R^M$]\label{inductive_limit}
Let $\{\R^M\}_{M \in \N}$ be the family of Euclidean spaces with $\N$ as their directed set. We call the direct limit of the sequence of spaces the inductive limit of convergence on $\R^M$. It will be denoted as $\Rinf$.
\end{definition}

\begin{definition}[Co-meagre/Residual]\label{Co_meagre}
Let $G$ be a topological space. We call a set $X$ meagre in $G$ if there exists a countable collection of nowhere dense subsets $N_n$ in $G$ s.t. $X = \cup_n ^\infty N_n$. 

We say $X$ is co-meagre or residual in $G$ if $X^c = G \setminus X$ is meagre. A property holding on a co-meagre/residual subset of $G$ is said to be generic.
\end{definition}

\begin{definition}[Fredholm operators and maps, Sec. 1 \cite{Smale_infinite_sard_65}]\label{Fredholm}
Let $T \in \cC^k(X, Y)$ be a bounded linear operator between Banach spaces. We say $T$ is Fredholm with index $k$ if: \begin{enumerate}
    \item $\dim(\ker(T)) < \infty$.
    \item Range of $T$ is closed.
    \item $\text{co}\ker (T)$ has finite dimension.
    \item $\dim(\ker(T)) - \dim(\text{co}\ker (T)) = k$.
\end{enumerate}
We say a (possibly nonlinear) mapping $F \in \cC^1(X, Y)$ is a Fredholm map of index $k$ if $DF: X \to Y$ is a Fredholm operator with index $k$ at each $x \in X$.
\end{definition}
\vspace{\baselineskip}

\subsection{Lojasiewicz Inequality and other key concepts for convergence analysis}
The following notions are key to considering convergence of optimisation flows:

\begin{definition}[Coercive mapping]\label{coercive}
    We will say a map $\CL: G \to \R^+$ over a metric space $G$ is coercive iff there exists an $\eps > 0$, s.t. the level set $\{g \in G: \CL[g] < \eps\}$ is a bounded subset of $G$.
\end{definition}
\begin{definition}[Norm Coercive mapping]\label{norm_coercive}
$\FF:G \to H$ is called a norm-coercive mapping if:
    \[
    \|g_i\|_G \to \infty \implies \|\FF[g_i]\|_H \to \infty
    \]

\end{definition}

We define the following gradations of convexity for a functional $\CL \in \cC^2(\CBep, \R)$:

\begin{definition}[Convex]\label{vanilla_con}
 $\CL \in \cC^2(\CBep, \R)$ is convex in $\CBep$, iff $\nabla^2\CL$ is a positive semi-definite operator (p.s.d.) for all $g \in \CBep$: i.e.\ $\sigma(\nabla^2\CL[g]) \subset \mbbR^+$.
\end{definition}

\begin{definition}[Strictly Convex]\label{strict_convex}
$\CL$ is strictly convex in $\CBep$, if for all $g \in \CBep - A$, $\sigma(\nabla^2\CL[g]) \subset \mbbR^+ \setminus \{0\}$, where $A$ is an at most countable set.
\end{definition}

\begin{definition}[$m-$Strongly Convex]\label{strong_convex}
  $\CL$ is $m-$strongly convex in $\CBep$ if $\sigma(\nabla^2\CL[g]) \subset [m, \infty]$ for all $g \in \CBep$. 
\end{definition}

\begin{definition}[Generalized Lojasiewicz Inequality]\label{General_LI}
    Let $\CL \in \cC^2(G)$, where $G$ is a separable Banach manifold. $\CL$ is said to satisfy the Lojasiewicz Inequality near $\Phi$ if there exists a neighbourhood $\CBep \subset G$ and $\alpha \in [0.5, 1), C > 0$ such that for all $g \in \CBep$:
    \begin{equation}\label{generalised_Banach_LI}
    |\CL[g] - \CL[\Phi]|^\alpha \leq C \|D\CL[g] \|_{(T_g G)^*} 
    \end{equation}
    
    Further, assume $G$ has sufficient structure to admit $T_g G \simeq (T_g G)^*$ under a $\cC^1$ and invertible duality or Riemannian metric $r_g: T_g G \times T_g G \to \R$, $r_g(\nabla\CL[g], v) = D\CL[g](v)$. We say that $\CL$ satisfies the gradient Lojasiewicz Inequality at $\Phi$ if there exists a neighbourhood $\CBep \subset G$ and $\alpha \in [0.5, 1), C > 0$ such that for all $g \in \CBep$:
    \begin{equation}\label{generalized_gradient_LI}    
    |\CL[g] - \CL[\Phi]|^\alpha \leq C \|\nabla \CL[g] \|_{T_g G} 
    \end{equation}
\end{definition}

\begin{remark}
    Note that if $G$ is a Hilbert manifold, Riesz representation generates a canonical identification between $D\CL[g]$ and $\nabla \CL[g]$, in turn allowing the kinds of canonical gradient flows given in Eq. \ref{idealized_G_flow}. However, if $G$ only admits a Banach–Riemannian or weaker structure, we need to pick the “right'' duality or $r_g$, since there is no canonical isomorphism we can use to transform Eq. \ref{generalised_Banach_LI} into \ref{generalized_gradient_LI}. $G$ will usually be a Hilbert space, unless otherwise mentioned.
\end{remark}

We will need the notions of semi-algebraic sets and o-minimality to show how the Lojasiewicz inequality is used to characterize convergence in certain applications:
\begin{definition}[Semi-algebraic Sets, Sec. 1, \cite{bier_milman_semianalytic_88}]\label{Semi-algebraic}
A subset $X \subset \R^k, k \in \N$ is called a semi-algebraic subset iff there exist polynomials $p_{ij}(s)$ and $q_{ij}(s)$, where $i = 1,...,p,$ and $j = 1, ...,q$ s.t.
\[
X = \bigcup_{i=1}^p \{s: p_{ij}(s) =0, q_{ij}(s) > 0, j = 1,...,q\}
\]
We will call a map semi-algebraic if its graph set is semi-algebraic set.
\end{definition}

\begin{definition}[Real analytic set]\label{real_analytic_set}
     $X \subset \R^k$ is called a real analytic set if it is the zero set of finitely many real analytic functions.
\end{definition}

\begin{lemma}[Curve Selection Lemma, Prop. 2.5.5 and Prop. 8.1.13 \cite{Curve_selection_Bochnak_98}]\label{curve_selection_lemma}
    Let $X \subset \R^k, k \in \N$ be a real-analytic set and $x \in \overline{X}$. Then, there exists a real analytic curve $\phi$ such that $\phi: [0, \epsilon) \to \R^k$, $\phi(0) = x$, and $\phi(t) \in X$ for all $t \in (0, \epsilon)$. Further, if $X \subset R^k, k \in \N$ is a semi-algebraic set, then a semi-algebraic curve $\phi$ exists as well.
\end{lemma}

\begin{definition}[o-minimality]\label{o_minimal_structure}
Let $S = \bigcup_{n\in \N} S_n$ where each $S_n$ is a family of subsets in $\R^n$. We say that $S$ is an o-minimal structure on $(R, +, \cdot)$ if:
\begin{enumerate}
    \item Each $S_n$ is a Boolean algebra (closed under finite set-theoretic operations)
    \item $A \in S_j, B \in S_k \implies A\times B \in S_{j+k}$
    \item If $A \in S_{j+k}$ and $\pi: R^{j+k}\to \R^j$ represents projection onto first $j$ coordinates, then $\pi(A) \in S_j$.
    \item If $p, q_1, ..., q_k \in \Q[X_1, ..., X_n]$, then $\{x \in \R^n: p(x) = 0, q_1(x) > 0, ... q_k(x) > 0 \} \in S_n$.
    \item $S_1$ is the collection of sets that are finite unions of open intervals and points.
\end{enumerate}
A set $X \in \R^n$ is called definable if $X \in S_n$. We will call a map definable in an o-minimal structure (or an o-minimal map in short), if its graph set is definable in an o-minimal structure.

In a very rough sense, we can associate a map being o-minimal as it satisfying the most general convergence condition variant of \ref{Lojasiewiczineqaulity} \cite{kurdyka98}.

\end{definition}

\subsection{Parameter Spaces and Architectures}\label{Rinfty}

The Euclidean/Hilbert spaces $\R^M, M \in \N$ are the natural domain spaces for $\CN$ with a finite number of parameters. The situation changes dramatically when we consider the notion of $M = \infty$. Ultimately, we need to define and rely on three separate types of sequence spaces: the vector space of all real sequences $\R^\N$, the Hilbert space of square summable real sequences $l^2$, and the LB space that is the inductive limit of convergence on $\R^M$, denoted as $\Rinf$.

$\Rinf$ as a set is also a linear subspace in $l^2$. Thus, we can also consider it using the subspace topology of $l^2$. We note that the closure of $\Rinf$ within the inductive limit topology is itself, while under the subspace topology, the closure is $l^2$. Working in the parametric space $\Rinf$ (with the inductive limit topology) raises issues since it is not a Banach space, metrizable, or norm-able. The classical Frechet derivative is thus not definable over such a space, but we may still speak of a $\cC^k$ map $\CN: \Rinf \to G$ in a precise sense. Formally:
\begin{definition}[\cite{KA97}]\label{C^k_Rinf}
    We say $\CN \in \cC^k(\Rinf, G)$ if for each $i \in \N$, we may write:
    \[
    \CN|_{\R^{M_i}} = \CN_i, \qquad D^j\CN|_{\R^{M_i}} = D^j\CN_i, \, \forall j \leq k
    \]
\end{definition}
\vspace{\baselineskip}

\subsection{Neural Tangent Kernels}
Since the Neural Tangent Kernels serve as an important inspiration and object for this work, we provide a deeper look at how they originated. Let us begin with objects originally termed Neural Tangent Kernel and contrast their heavy reliance on the data and input as a control for optimisation dynamics.

\subsubsection{Conventional Neural Tangent Kernels}\label{NTK_intro}
The Neural Tangent Kernel (NTK) was introduced in \cite{jacot_ntk_18} and has been shown to serve as a control for optimization dynamics in some idealized settings. There are several related objects that the literature may be alluding to when using the term “NTK", depending on how the parameters of the architecture are initialised, the relationship to the sample set, etc. Here, we note the two most common forms (in the main body, we work with our interpretation).

A standard definition views the NTK as the kernel for the feature map $(\cdot)\to \frac{\partial}{\partial w_i} \cNw (\cdot)$. Formally, let $\CN(\weightvector)(\cdot): \R^M \times \R^m \to \R^n$ represent the model as a map between $m, n$ dimensional Euclidean spaces. Then, the NTK is the following object, in terms of some $x, y \in \R^m$:
\begin{equation}\label{NTK_og_def}
    \Theta_{n_1, n_2}(\weightvector)(x, y) = \sum_1^M \Bigg[\frac{\partial}{\partial {w_i}}\CN_{n_1}(\weightvector)(x)\Bigg]\Bigg[\frac{\partial}{\partial {w_i}}\CN_{n_2}(\weightvector)(y)\Bigg]
\end{equation}
where $\CN_{n_i}(\weightvector)(x)$ represents the $n_i^{th}$ output of the model for an input $x$. This object was shown to be a control for optimization dynamics under the class of empirical risk minimization problems.

Another essentially equivalent formulation, from \cite{ntk19jae}, involves \textit{unrolling} or concatenating the model function over the set of sample inputs $\CX$, so that we have a function of models over each $x_a \in \CX$. We then create a Jacobian of this vector of functions with respect to the parameters. Assuming our model is a $J$ output map from the domain of interest $\tmcX$ to $\tilde{\CY}$: \[
\CN^\CX = 
\begin{pmatrix}
\CN^{x_1}_1 \\
: \\
\CN^{x_N}_1 \\
\CN^{x_1}_2 \\
: \\
\CN^{x_N}_J
\end{pmatrix},
\qquad \qquad \qquad
\NdwX = 
\begin{pmatrix}
\frac{\partial \CN^{x_1}_1}{\partial w_1} & ... & \frac{\partial \CN^{x_1}_1}{\partial w_M} \\
: & ... & :\\
\frac{\partial \CN^{x_N}_1}{\partial w_1} & ... & \frac{\partial \CN^{x_N}_1}{\partial w_M} \\
\frac{\partial \CN^{x_1}_2}{\partial w_1} & ... & \frac{\partial \CN^{x_1}_2}{\partial w_M} \\
: & ... & :\\
\frac{\partial \CN^{x_N}_J}{\partial w_1} & ... & \frac{\partial \CN^{x_N}_J}{\partial w_M} 
\end{pmatrix},
\qquad\qquad\qquad
\Theta^\CX = \NdwX(\CN^{\CX}_{\weightvector})^\dagger
\]

In this formulation, the NTK is a $J N\times J N$ symmetric matrix.

\subsection{Stochastic Gradient Flows}\label{SGD_uses_lin_vs_nlin_N}
A well-studied and used stochastic gradient flow that is of relevance to Eq. \ref{param_GD} is given by:

\begin{equation}\label{param_SGDyn}
\begin{aligned}
& \frac{d }{dt}\wt = -\nabla \cL[\wt] + \alpha(t)\beta(t)\xi(t), \qquad \qquad \alpha(t) = \sqrt{\frac{c}{\log(2+t)}}\\
\implies 
& {\frac{d}{dt}\CNt = -\Theta(t)\nabla\CL[\CNt]}  + \alpha(t)\Bigg[\frac{\beta^2(t)}{2}\CNww + \Ndw\beta(t)\xi(t)\Bigg]\\
\end{aligned}
\end{equation}
where $\xi(t)$ is a Gaussian white noise term, $\alpha(t)$ is the annealing schedule, and $\beta(t)$ is the standard deviation. The corresponding Fokker-Planck equation for the density $\rho_\weightvector(t)$ is:
\begin{equation}\label{FPE_SGD}
\frac{d}{dt}\rho_\weightvector(t) = \nabla \cdot \Bigg(\rho_\weightvector(t) \DcLw  + \alpha(t)\frac{\beta^2(\weightvector, t)}{2} \nabla \rho_\weightvector(t)  \Bigg)    
\end{equation}

Under Eq. \ref{param_SGDyn}, if $\xi$ is a Brownian motion with constant $\beta$, $\CL$ is coercive and well-behaved, and $\CBMep$ is non-empty, almost every $\wzero$ can lead to a well-initialised $\cNw$: we have a convergence in law to the global minimizer of $\cL$ \cite{sgd_global_convg_geman_86}. Once we are near that global minimizer, we can drop the noise terms to recover Eq. \ref{param_GD} with a well-initialised model (we can view the annealing schedule $\alpha(t)$ as the term doing this for us).

Unfortunately, this convergence is not almost surely (it does not happen with probability 1 for all initial conditions). Further, logarithmic annealing schedules are costly and inverse polynomial ones are used in practice: it is common practice to obtain coarse estimates and then find better models through finer methods \cite{good_guess_then_refine_Ilya_13, good_guess_k_Means_13, warm_start_Yildrim_02}.

\end{document}